\documentclass[11pt]{article}
\usepackage{amssymb}
\usepackage{latexsym}

\newtheorem{thm}{Theorem}[section]
\newtheorem{lem}[thm]{Lemma}
\newtheorem{cor}[thm]{Corollary}
\newtheorem{pro}[thm]{Proposition}
\newtheorem{ex}[thm]{Example}
\newtheorem{rmk}[thm]{Remark}
\newtheorem{defi}[thm]{Definition}

\newcommand{\gm }{\Gamma }
\newcommand{\lon }{\longrightarrow }
\newcommand{\be }{\begin{eqnarray*}}
\newcommand{\ee }{\end{eqnarray*}}

\newcommand{\poidd }[2]{#1\gpd #2}
\newcommand{\poiddd }[3]{ (#1\gpd #2, \alpha_{#3}, \beta_{#3})}

\newcommand{\pf}{\noindent{\bf Proof.}\ }
\newcommand{\qed}{\begin{flushright} $\Box$\ \ \ \ \ \
                    \end{flushright}}
\newcommand{\complex}{{\mathbb C}}
\newcommand{\reals}{{\mathbb R}}

\newcommand{\frakg}{{\mathfrak g}}
\newcommand{\frakh}{{\mathfrak h}}
\newcommand{\frakl}{{\mathfrak l}}
\newcommand{\frakm}{{\mathfrak m}}

\newcommand{\frakk}{{\mathfrak k}}
\newcommand{\frakt}{{\mathfrak t}}

\newcommand{\half}{\frac{1}{2}}

\newcommand{\pr}{\mbox{pr}}

\newcommand{\calx}{{\mathfrak X}}

\newcommand{\smalcirc}{\mbox{\tiny{$\circ $}}}
\newcommand{\ra}{\rangle}
\newcommand{\la}{\langle}   

\def\description label#1{\hfil\bf[#1]\hfil}
\newcommand{\parrr}[2]{\frac{\partial #1}{\partial #2}}

\newcommand{\Q}{Dirac }
\newcommand{\sy}{symmetric }

\newcommand{\ii}{\sqrt{-1}}

\def\sdp{\mathbin{\hbox{$\mapstochar\kern-.3333em\times$}}}
\def\pds{\mathbin{\hbox{$\times\kern-.55em\mapstochar\,$}}}

\newcommand{\wed}{\mathbin{\lower1.5pt\hbox{$\scriptstyle{\wedge}$}}}

\let\Tilde=\widetilde

\let\Vec=\overrightarrow
\let\ceV=\overleftarrow

\def\chigh{{\raise1.5pt\hbox{$\chi$}}}
\let\phi=\varphi
\def\til0{\Tilde{0}}

\def\dminus{\raise2pt\hbox{\vrule height1pt width 2ex}\hskip3pt}

\def\pback#1{\mathbin{{{\lower1.2ex\hbox{$\times$}}\atop #1}}}

\def\vlra{\hbox{$\,-\!\!\!-\!\!\!-\!\!\!-\!\!\!-\!\!\!
-\!\!\!-\!\!\!-\!\!\!-\!\!\!-\!\!\!\longrightarrow\,$}}

\def\gpd{\,\lower1pt\hbox{$\longrightarrow$}\hskip-.24in\raise2pt
               \hbox{$\longrightarrow$}\,}

\def\lgpd{\,\lower1pt\hbox{$\vlra$}\hskip-1.02in\raise2pt\hbox{$\vlra$}\,}

\def\llgpd{\,\lower1pt\hbox{$\vvlra$}\hskip-1.3in\raise2pt\hbox{$\vvlra$}\,
}

\begin{document}

\title{{\bf \Q submanifolds and Poisson involutions}}
\author{ PING XU \thanks{ Research partially supported by   NSF
          grant DMS00-72171.}\\
   Department of Mathematics\\The  Pennsylvania State University \\
University Park, PA 16802, USA\\
          {\sf email: ping@math.psu.edu } }

\date{}

\maketitle
\centerline{Dedicated to Rencontres   Math\'ematiques de Glanon
on the occasion of her  fifth birthday}

\begin{abstract}
Dirac submanifolds are a natural generalization in the Poisson category
 for symplectic submanifolds of a symplectic manifold.  In
a certain sense  they correspond to symplectic subgroupoids
 of the symplectic groupoid  of the given Poisson manifold.
 In particular, Dirac  submanifolds arise as the stable locus of a  Poisson involution.
In  this paper, we provide a general study for these
submanifolds including both local and global aspects.

In the second part of the paper,
we  study Poisson involutions and the induced Poisson
structures on their stable locuses. We discuss the
Poisson involutions on a special class of  Poisson groups, and more
generally  Poisson groupoids, called symmetric Poisson
groups (and symmetric Poisson groupoids). Many well-known examples,
including the  standard Poisson group structures on semi-simple Lie groups, Bruhat
Poisson structures on compact semi-simple Lie groups, and
Poisson groupoids connecting with dynamical $r$-matrices of semi-simple
Lie algebras are symmetric, so they admit a Poisson involution.
 For symmetric Poisson  groups, the relation between
 the  stable locus Poisson structure
and Poisson symmetric spaces is  discussed. As a consequence,
we show that the  Dubrovin-Ugaglia-Boalch-Bondal Poisson
structure on the space of Stokes matrices $U_{+}$
appearing in Dubrovin's theory of Frobenius manifolds is  indeed a Poisson
symmetric space for the Poisson group $B_{+}*B_{-}$.
\end{abstract}

\section{Introduction}

The underlying structure of any Hamiltonian system is a Poisson manifold.
To deal with mechanics with constraints, it is always desirable to understand
how to put a Poisson structure on a submanifold of a Poisson manifold.
 A naive way is  to consider
Poisson submanifolds. However, these are not too much different from
the original Poisson manifold from the viewpoint of the Hamiltonian
systems. On the other hand, for symplectic manifolds, there
do not exist any nontrivial Poisson submanifolds.
However Dirac was able to write down a Poisson bracket
for  a submanifold of a symplectic manifold which is
given  by  a set of constraints:
\begin{equation}
Q=\{x\in P|\phi_{i}(x)=0, i =1, \cdots , k\}
\end{equation}
such that the matrix $(\{\phi_{i}, \phi_{j}\})$ is invertible on $Q$.
This is the famous Dirac bracket \cite{Dirac}. In this case, $Q$ is
 a symplectic submanifold, i.e., the pull back on $Q$ of the symplectic
form is non-degenerate.

There has appeared a lot of work attempting to generalize Dirac brackets,
 for example, the notion of cosymplectic manifolds of Weinstein
\cite{W83}, Poisson reduction of Marsden-Ratiu \cite{MR}, just to name a few.
In particular, Courant presented a unified approach  to this question
 by introducing the notion of Dirac structures \cite{Courant}, by
which one could obtain a Poisson bracket on admissible
functions on a submanifold $Q$. In some situation, one
indeed can get a Poisson structure on all functions on $Q$.
Then $Q$ becomes a Poisson manifold itself.

In his  study of  Frobenius manifolds, which is connected with
2-dimensional topological quantum field theories, Dubrovin  recently
found a Poisson structure on  $U_{+}$, 
the space of  upper triangular matrices with ones on the diagonal,
by viewing it as a space of
Stokes matrices. Indeed  Dubrovin  identifies $U_{+}$  with the
local moduli spaces of semisimple Frobenius manifolds.
 In particular, an explicit formula
was found for the Poisson bracket in the three dimensional case \cite{Dubrovin}:
\begin{equation}
\pmatrix{
1&  x& y \cr
0 & 1 &z\cr
0& 0& 1}     
\end{equation}
\begin{equation} 
\{x, y\}=xy-2z, \ \{y, z\}=yz-2x, \ \{z, x\}=zx-2y.
\end{equation} 
Such a Poisson structure admits various nice properties. For instance,
it naturally  admits a  braid group action.
The casimir function is the Markoff polynomial $x^2 +y^2 +z^2 -xyz$.
Its linear and quadratic parts give rise to a biHamiltonian
structure, etc.
Then Ugaglia extended Dubrovin's formula to the $n\times n$ 
case \cite{Ugaglia}.
Recently, in connection with his study of  the monodromy  map,
Boalch \cite{Boalch} proved that $U_+$ arises as
the stable locus of a Poisson involution on the Poisson
group $B_{+}*B_{-}$ and that the above Poisson structure on
$U_+$ is induced from the standard Poisson structure on $B_{+}*B_{-}$.

>From a completely different angle,  independently Bondal  
discovered exactly the same Poisson structure on $U_{+}$
 in his study of  the theory of triangulated 
categories \cite{Bondal1}. He also studied extensively this
 Poisson structure including the braid group action and
symplectic leaves etc. In his approach, instead of writing down the
Poisson structure on $U_{+}$,  first of all Bondal discovered
a symplectic groupoid ${\cal M}$ whose space of objects is $U_{+}$. Then
the general theory of symplectic groupoids \cite{W87} implies
that $U_{+}$ is a Poisson manifold.  What is more interesting is,
in a sebsequent paper \cite{Bondal2},
he discovered  an extremely
  simple connection between his symplectic groupoid  ${\cal M}$
and the standard symplectic groupoid  $\gm$ over the  Poisson group $B_{+}*B_{-}$
of  Lu-Weinstein  \cite{LW1}. Namely, ${\cal M}$ is simply
a symplectic subgroupoid of $\gm$ which can be realized as
the stable locus of an involutive symplectic  groupoid automorphism of
$\gm$.

Bondal's work suggests a simple   fact, which  was  somehow overlooked in
the literature, namely a submanifold inherits a  natural Poisson structure if
it can be realized as the base space of a symplectic subgroupoid.
A natural question arises as to what are these submanifolds and
how they can be characterized. One of the  main purpose of the paper
is to answer this question. These submanifolds will be called
 \Q submanifolds. Symplectic subgroupoids are very simple to
describe: they are subgroupoids and in the mean time
symplectic submanifolds.  In contrast, Dirac submanifolds   
are not so simple as we shall see.
 There are some interesting and
rich geometry there (both global and local), which we 
believe  deserve
further studies.

\Q submanifolds are a special case  of those submanifolds,
 on which  the admissible functions for the
pulled back Dirac structure happen to be all functions
in terms of \cite{Courant, MR}. In other words,
the intersection of a \Q submanifold
with symplectic leaves of $P$ are symplectic submanifolds of the leaves.
This feature explains where the induced Poisson structure comes
from for a \Q submanifold. However, not all such submanifolds
are \Q submanifolds. For instance, symplectic leavess (except for the 
zero point) of ${\mathfrak su}(2)$ are not  \Q submanifolds.
It is still not clear at the moment how to describe
the global obstruction  in general.
On the other hand, when the Poisson manifold is symplectic, 
\Q submanifolds  are precisely
symplectic submanifold.  Other examples include
cosymplectic submanifolds   and stable locus of a Poisson involution.

The second aim  of the paper is to study systematically 
Poisson involutions and the induced Poisson structures
on stable locuses.
When the underlying Poisson manifolds are Poisson  groups or
more generally Poisson groupoids, there is an effective way
of producing a Poisson involution, namely through
their infinitesimal invariants: Lie bialgebras or Lie
bialgebroids. These are called \sy Poisson groupoids
and \sy  Lie bialgebroids. As we see, such a Poisson
involution exists in almost every well-known
example of Poisson groupoids and Poisson groups,
including the standard Poisson group structures on semi-simple Lie groups, Bruhat
Poisson structures on compact semi-simple Lie groups, and
Poisson groupoids connecting with dynamical $r$-matrices of semi-simple
Lie algebras.                    
  For Poisson groups,  they were  studied by   Fernades 
in connection with  Poisson symmetric spaces \cite{Fern, Fern1},
i.e., symmetric spaces which are Poisson homogeneous spaces.
It turns out that the induced Poisson structure on the
stable locus $Q$  of the  Poisson involution of a \sy Poisson group
is closely connected with Poisson symmetric spaces.
In particular, we prove that the identity connected 
component of $Q$ is always  a Poisson symmetric space.
As a consequence, we show  that the DUBB-Poisson structure
on Stokes matrices $U_{+}$ is a Poisson symmetric space for
the Poisson group $B_{+}*B_{-}$.

The paper is organized as follows.  In Section 2, 
we introduce the definition of \Q submanifolds and
study their basic properties. Local \Q submanifolds
are also introduced and their connection with transverse
Poisson structures is discussed. Section 3 is devoted to
the study of some further properties. In particular, we
study how the  modular class of a \Q submanifold is related
to that of the Poisson manifold $P$.  We also study Poisson
group actions on  \Q submanifolds. Finally we prove that
\Q submanifolds  are indeed infinitesimal version of 
symplectic subgroupoids. In Section 4, we investigate 
stable locuses of  Poisson involutions,  and study
Poisson involutions on Poisson groupoids by introducing
the notion of \sy Poisson groupoids. In Section 5,
we consider particularly \sy   Poisson groups and
the induced Poisson structures on stable locuses.
The connection with Poisson symmetric spaces is 
discussed.

We remark that one should not confuse \Q submanifolds here with
the notion of \Q manifolds of Courant \cite{Courant}.
Courant's \Q manifolds are manifolds equipped with a
Dirac structure, which generalize the notion of  both Poisson and
presymplectic manifolds. In an earlier version, some
other  names such as $Q$-submanifolds and IR-submanifolds
were suggested, but we feel that neither of these names
 reflects the complete nature of the objects we study here.
At the end, we decided to call them \Q submanifolds, which
at least contains a famous name that people have heard of.

{\bf Acknowledgments.}
The author would like to thank ESI and University of Pennsylvania
for their hospitality while work on this project was being done.
He also wishes to thank Philip Boalch, Alexei Bondal, Sam Evens,
 Rui Fernades, Friedrich  Knop,  Jiang-hua Lu and Alan Weinstein for  numerous 
fruitful discussions, comments  and suggestions. In particular,
he is  grateful to Alexei Bondal for allowing him to be accessible
to his unpublished manuscripts \cite{Bondal1, Bondal2}.
Finally special thanks go to all  his friends in Glanon,
in particular, Fr\'ed\'eric Bidegain,  S\'ebastien Mich\'ea,
and  Fran\c cois Nadaud, 
from whom he  benefited a lot not only mathematically, but also 
learned much more beyond.

\section{\Q submanifolds}

This section is devoted to the study of general aspects of \Q submanifolds.

\subsection{Definition and properties}

Let us introduce the definition first.

\begin{defi}
\label{def:Q}
A submanifold $Q$ of a Poisson manifold  $P$ is called
a \Q submanifold if the tangent bundle of $P$ along $Q$  admits
a  vector bundle decomposition: 
\begin{equation}
\label{eq:decom}
T_{Q}P=TQ\oplus V_{Q}
\end{equation}
 so that  $V_{Q}^{\perp}$ is a Lie subalgebroid of $T^* P$,
where   $T^*P$  is equipped with the standard cotangent
bundle Lie algebroid structure.
\end{defi}

Note that  the last condition above is equivalent to
 that $V_{Q}\subset TP$ is a coisotropic submanifold of the tangent Poisson
manifold $TP$.  So alternatively,  we have

\begin{pro}
A submanifold $Q\subset P$ is a \Q submanifold
iff there is a decomposition as in Equation (\ref{eq:decom}) so that
 $V_{Q}\subset TP$ is a coisotropic submanifold of the tangent Poisson
manifold $TP$.
\end{pro}

In what follows, we will see that $Q$ itself must be a Poisson manifold.
 However, $Q$  in general is not a Poisson submanifold.
We need to introduce some notations. By $\pr$, we denote the
bundle map $T_{Q}P\lon TQ$ obtained simply  by taking the
projection along the decomposition (\ref{eq:decom}).
And let   $\pr^* : T^*Q\lon T^*P$ denote 
the dual  of $\pr$ by considering $T_{Q}^* P$ as a subbundle of
$T^*P$.  By $\pr_{*}$, we denote the map from $\calx^{d}(P)$
to $\calx^{d}(Q)$ naturally  induced from $\pr$, which is defined
by $\pr_{*} (D)=\pr ( D|_{Q})$, $\forall D\in \calx^{d}(P)$.

We summarize  some important properties of \Q submanifolds in the following

\begin{thm}
\label{thm:Q}
Let $Q$ be  a $\Q$ submanifold of a Poisson manifold $(P, \pi )$.
Then
\begin{enumerate}
\item $\pi|_{Q}=\pi_{Q}+\pi'$, where $\pi_{Q} \in \gm ( \wedge^{2}TQ)$ and
$\pi' \in \gm (\wedge^{2}V_{Q})$;
\item  $\pi_{Q}$ is a Poisson tensor on $Q$;
\item  $\pr^{*}: T^* Q\lon T^* P$ is a
Lie algebroid morphism, where both $T^* Q$ and $T^*P$
are equipped with the cotangent bundle Lie algebroid structures;
\item for any $X \in \calx (P)$,
\begin{equation}
\label{eq:pr}
 \pr_*  [X , \pi ]=[ \pr_* X , \pi_{Q}];
\end{equation}
\item for any $x\in Q$, 
 $\pi_{Q}^{\#} (T_{x}^{*}Q)=\pi^{\#}(T_{x}^*P) \cap T_{x}Q$;
\item  for any $x\in Q$,  $\pi_{Q}^{\#} (T_{x}^{*}Q)$
is a  symplectic subspace of $\pi^{\#}(T_{x}^*P)$. Hence,
  each symplectic leaf of $Q$ is the  intersection of $Q$ with
a  symplectic leaf of $P$, which is  a symplectic submanifold
of that leaf.
\end{enumerate}
\end{thm}          

Before proving this theorem, we need a couple of lemmas.
The following lemma, which  can
also  be easily verified directly, follows  from the fact that the natural
inclusion $TQ\lon TP$ is a Lie algebroid morphism.

\begin{lem}
\label{lem:bracket}
Let $Q\subseteq P$ be  a submanifold. Assume that $D\in \calx^{d}(P)$
and  $D'\in \calx^{d'}(P)$ are multi-vector fields 
 tangent to $Q$, i.e., $D|_{Q}\in \calx^{d}(Q)$ and
 $D'|_{Q}\in \calx^{d'}(Q)$. Then
$$[D, \ D']|_{Q}=[D|_{Q}, \ D'|_{Q}].$$
\end{lem}

Next is the following

\begin{lem}
\label{lem:piQ}
Assume that $Q$ is a submanifold of a Poisson manifold $(P, \pi )$
such that there is a vector bundle decomposition
 $T_{Q}P=TQ\oplus V_{Q}$.  Moreover assume that
 $\pi |_{Q}=\pi_{Q}+\pi '$, where
$\pi_{Q}\in \gm (\wedge^{2}TQ )$ and $\pi '\in \gm (\wedge^{2} V_{Q})$.
Then $\pi_{Q}$ is a Poisson tensor on $Q$.
\end{lem}
\pf Write $\pi ' =\sum_{i} X_{i}\wedge Y_{i}$ where $X_{i}, Y_{i}\in
\gm (V_{Q} )$. Now let $\tilde{X}_{i}, \tilde{Y}_{i}\in
\calx (P)$ be   (local) extensions  of  $X_{i}, Y_{i}$, and
$\tilde{\pi }'=  \sum_{i} \tilde{X}_{i}\wedge \tilde{Y}_{i}$.
Let $\pi ''=\pi -\tilde{\pi }'\in \calx^{2}(P)$. Then clearly
$\tilde{\pi }'|_{Q}=\pi '$ and  $\pi ''|_{Q}=\pi_{Q}$.
Now it follows from $[\pi, \pi ]=0$ that
$[\pi '' , \pi '' ]=-2[\pi '' , \tilde{\pi }']-[\tilde{\pi }', \tilde{\pi }']$.
On the other hand, it is clear by definition that
$\pr_{*} [\pi '' , \tilde{\pi }']=\pr_{*} [\tilde{\pi }', \tilde{\pi }']=0$.
Thus $\pr_{*} [\pi '' , \pi '' ] =0$.
According to Lemma \ref{lem:bracket},   the latter implies
 that $[\pi_{Q} , \pi_{Q} ]=\pr_{*} [\pi '' , \pi '' ]=0$. 
This concludes the proof. \qed
{\bf Proof of Theorem \ref{thm:Q}} 
  By definition, $V_{Q}^{\perp}$ is a Lie subalgebroid of the
cotangent Lie algrebroid $T^*P$. By identifying $T^*Q$ with
 $V_{Q}^{\perp}$, one obtains a Lie algebroid structure on
$T^*Q$, and a Lie algebroid morphism $\phi: T^*Q\lon T^*P$. Clearly,
$\phi=\pr^*$. By  $\rho_{Q}$, we denote the anchor map 
of the Lie algebroid $T^*Q$. Thus we have $i \smalcirc \rho_{Q}=
\pi^{\#}\smalcirc \phi$, where $i: TQ\lon TP$ is the natural inclusion. 
It  follows that $\rho_{Q}=\pr \smalcirc i\smalcirc \rho_{Q}=
\pr \smalcirc \pi^{\# }\smalcirc \phi
=\pr \smalcirc \pi^{\# }\smalcirc \pr^{*}$. Hence $\rho_{Q}: T^*Q \lon
TQ $ is skew-symmetric. Thus it  defines a bivector field
$\pi_{Q}\in \gm (\wedge^{2}TQ)$ so that $\rho_{Q}=
\pi_{Q}^{\#}$. Under the decomposition (\ref{eq:decom}),
 we have $\wedge^2 T_{Q}P=\wedge^{2}TQ\oplus (TQ\wedge  V_{Q})\oplus
\wedge^{2}V_{Q}$. It is clear that $\pi_{Q}$ is the
$\gm (\wedge^2 TQ )$-part of $\pi|_{Q}$ under the above decomposition.
Since $\pi^{\#}(V_{Q}^{\perp})\subset TQ$,  $\pi|_{Q}$ does not involve
any  mixed term, i.e., the $\gm (TQ\wedge  V_{Q})$-part.
Hence we have     $\pi|_{Q}=\pi_{Q}+\pi'$ with $\pi'
\in \gm (\wedge^{2}V_{Q})$.  This proves (i).

By Lemma \ref{lem:piQ}, $\pi_{Q}$ is indeed a Poisson tensor 
on $Q$.  Hence (ii) follows.
 Next we need to show that the  Lie algebroid   structure
on $T^*Q$ is indeed  the cotangent Lie algebroid
corresponding to
 the Poisson structure $\pi_{Q}$. Since $\phi$ is a Lie algebroid morphism,
 $\phi^*=\pr$ induces a morphism of the (graded) differential
 algebras $\pr_{*}: (\gm (\wedge ^* TP), d_{*P})\lon  
 (\gm (\wedge ^* TQ), d_{*Q})$. Since $T^* P$ is the cotangent 
Lie algebroid,  we know that
 $d_{*P}=[\pi,  \cdot ]$.  To prove the claim, it suffices to show that
$d_{*Q}=[\pi_Q ,  \cdot ]$. To this end,
 given any $X\in \calx (Q)$, choose an extension  $\tilde{X}
\in \calx(P) $.
Write $\pi = \pi'' + \tilde{\pi}'$
as in Lemma \ref{lem:piQ} so that $\pi'' |_{Q}=\pi_{Q}$ and
$\tilde{\pi}'|_{Q}\in \gm (\wedge^2 V_Q )$. Then 
$d_{*Q}X =(d_{*Q}\pr_{*})\tilde{X}=(\pr_{*} d_{*P})\tilde{X}
 =\pr_*  [\pi , \tilde{X}]=
\pr_*  [\pi'' + \tilde{\pi}', \tilde{X}]
=\pr_* [\pi''  , \tilde{X}] =[\pi_{Q}, X]$, where the last step
follows from Lemma \ref{lem:bracket}.
 This proves (iii), and therefore
(iv) as a consequence.

Next we  prove the relation
  $\pi_{Q}^{\#} (T_{x}^{*}Q)=\pi^{\#}(T_{x}^*P) \cap T_{x }Q$.
Since  $i \smalcirc \pi_{Q}^{\#}= \pi^{\#}\smalcirc \phi$, it is 
obvious that $ \pi_{Q}^{\#} (T_{x}^{*}Q) \subset
 \pi^{\#}(T_{x}^*P) \cap T_{x }Q$. 
Conversely, let $v\in  \pi^{\#}(T_{x}^*P) \cap T_{x }Q$ be
any vector. Then $v=\pi^{\#} \xi $ for some  $\xi \in T^{*}_{x}P$.
Now since  $T_{x}^{*}P=T_{x}Q^{\perp} \oplus V_{x}^{\perp}$, we may write
$\xi =\xi_{1}+\xi_2$, where $\xi_{1}\in T_{x} Q^{\perp}$ and
$\xi_2 \in V_{x}^{\perp}$. Since $V_{x}^{\perp} =\phi (T_{x}^{*}Q)$,
$\pi^{\#} \xi_{2}\in (\pi^{\#}\smalcirc \phi ) (T_{x}^{*}Q)=\pi_{Q}^{\#} 
(T_{x}^{*}Q)\subset  T_{x}Q$. Hence $ \pi^{\#} \xi_{1} =v-\pi^{\#} \xi_{2}\in 
T_{x}Q$. On the other hand, it is clear that $\pi^{\#} \xi_{1}
\in V_{x}$. Hence, $\pi^{\#} \xi_{1}=0$ and therefore $v=\pi^{\#} \xi_{2}
\in \pi_{Q}^{\#} (T_{x}^* Q)$.  Thus we have  proved
the relation
  $\pi_{Q}^{\#} (T_{x}^{*}Q)=\pi^{\#}(T_{x}^*P) \cap T_{x }Q$,
which implies that the symplectic leaves of $Q$ are
the intersection of the symplectic leaves 
of $P$ with $Q$. 

Finally,  let $D_{x}=\pi_{Q}^{\#} (T_{x}^{*}Q)$ and
$D_{x}'=\pi'^{\#} (T_{x}^{\perp}Q)$. It is simple to see
that $\pi^{\#} (T_{x}^*P)=D_{x}\oplus D_{x}'$, and 
$\pi_{Q}(x)\in \wedge^2 D_{x}$ and $\pi' (x)\in \wedge^2 D_{x}'$ are both
nondegenerate. Thus $\pi_{Q}(x)|_{D_{x}}^{-1} \oplus \pi'(x)|_{ D_{x}'}^{-1}$
is the inverse of $\pi (x)$ when being restricted
to $\pi^{\#} (T_{x}^*P)$. It follows that
$\pi_{Q}^{\#} (T_{x}^{*}Q)$ is indeed 
 a  symplectic subspace of $\pi^{\#}(T_{x}^*P)$. 
 This implies that any symplectic leaf of $Q$
is indeed a symplectic submanifold of a 
symplectic leaf of  $P$.
This concludes our proof of the theorem. \qed

As an immediate consequence, we have

\begin{cor}
\label{cor:2.6}
Assume that $Q$  is  a \Q submanifold of a Poisson manifold $P$.
Then we have
\begin{enumerate}
\item there is a morphism  on the level of Poisson cohomology
$$\pr_{*} : H^{*}_{\pi} (P)\lon  H^{*}_{\pi_Q} (Q);$$
\item if $X\in \calx (P)$ is a vector field such that
$X|_{Q}\in \gm (V_{Q})$, then  $\pr_{*} [X, \pi ]=0$.
\end{enumerate}
\end{cor}

\begin{rmk}
{\em From Theorem \ref{thm:Q}  (vi), we see that the choice
of the complementary $V_Q$ is immaterial for the purpose
of getting the Poisson structure on $Q$. Indeed, any
submanifold  whose intersections with symplectic
leaves of $P$ are symplectic submanifolds of
the leaves admits
a potential Poisson tensor, which, however, might be discontinuous.
This is simply the bivector field obtained by taking the
inverse of the  restriction of the leafwise
 symplectic   form to $Q$. In terms of the language of Dirac
structures, such  submanifolds precisely correspond  to those
for which  the pulled back Dirac structure \cite{Courant} of the
one corresponding to the graph of the Poisson tensor   on $P$ is
a bivector on each tangent space.  
In general, it might be discontinuous though. However, note that
even when it is smooth so that one obtains a 
Poisson structure on $Q$, it may still not be a \Q submanifold.
See Example \ref{ex:2.16} below. 

We also note that \Q submanifolds are a special case of 
the situation in \cite{MR}, where  general Poisson reduction
was studied. This provides another route to obatin the
Poisson structures on these submanifolds.}
\end{rmk}

Next proposition gives   an alternative definition
of \Q submanifolds, which is  presumably easier
to check in practice.
 
\begin{pro}
\label{pro:sufficient}
A submanifold $Q$ of a Poisson manifold $(P, \pi )$ is a \Q submanifold
if the following conditions are all satisfied:
\begin{enumerate}
\item  there is a  vector bundle decomposition $T_{Q}P=TQ\oplus V_{Q}$;
\item  $\pi |_{Q}=\pi_{Q}+\pi '$, where
$\pi_{Q}\in \gm (\wedge^{2}TQ )$ and $\pi '\in \gm (\wedge^{2} V_{Q})$.
\item for any $X' \in \gm (V_{Q})$, there is an extension
$X\in \calx (P)$ of $X'$  such that $\pr_*  [X , \pi ] =0$.
\end{enumerate}
\end{pro}
\pf From (i)-(ii), we know that $\pi_{Q}$ is
a Poisson tensor on $Q$, and therefore $T^* Q$
is a Lie algebroid.
The decomposition (i) induces  a natural identification
between $V_{Q}^{\perp}$ and $T^* Q$, which equips $V_{Q}^{\perp}$
with a Lie algebroid structure by pulling back  the cotangent Lie algebroid
on $T^* Q$. It remains to show that
this Lie algebroid structure on  $V_{Q}^{\perp}$
is indeed a Lie subalgebroid of $T^*P$. To this end,
it suffices to prove Equation (\ref{eq:pr}) for
any vector field $X\in \calx (P)$.

If $X\in \calx (P)$ such that $X|_{Q}$
is tangent to $Q$, Equation (\ref{eq:pr}) follows from
Lemma \ref{lem:bracket}.  On the
other hand, assume that  $X|_{Q}\in \gm (V_{Q})$. Then
$\pr_* [X , \pi ]= \pr_*  [X, \pi'' +\tilde{\pi} '] =\pr_*  [X, \pi'' ]$,
where $\pi''$ and $\tilde{\pi} '$ are  the
bivector fields as in the proof of Lemma \ref{lem:piQ}.
Since $\pi'' |_{Q}
=\pi_{Q}$ is tangent to $Q$, $ [X, \pi'' ]|_{Q}$  depends only
on $X|_{Q}$. From assumption (iii), we thus have $\pr_* [X , \pi ] =0$.
This concludes the proof.\qed

\begin{rmk}
\label{rmk:sufficient}
{\em Conditions (ii) and (iii) in Proposition \ref{pro:sufficient}
 can be replaced, respectively,  by the following equivalent conditions:

(ii'). $\pi^{\#} (V_{Q}^{\perp})\subseteq TQ$;

(iii'). for any point $x$ in $Q$, there are a set of local vector fields
$X_{1}, \cdots , X_{k}\in \calx (P)$ around $x$ such that 
$X_i|_{Q} \in \gm (V_{Q} )$,
$i=1, , \cdots , k$  consists of a fiberwise basis for $V_{Q}$
and satisfies the property that $\pr_* [X_i , \pi ] =0$, $i=1, , \cdots , k$.}
\end{rmk}

Recall that cosymplectic submanifolds of a Poisson manifold
$P$  are those,
which  are characterized by the two properties  \cite{W83}:
\begin{enumerate}
\item $Q$ intersects each symplectic leaf of $P$  transversely;
\item at each point of $Q$, the intersection of $TQ$ with
the tangent space of the symplectic leaf is a symplectic
subspace.
\end{enumerate}                  

\begin{lem}
\label{lem:co}
A  submanifold $Q$ of a Poisson manifold $(P, \pi )$ is
cosymplectic iff it satisfies the conditions (i)-(ii) as in 
Proposition \ref{pro:sufficient} with the property that
 $\pi ' \in \gm (\wedge^{2} V_{Q})$ is non-degenerate. 
\end{lem}
\pf If $Q$ is cosymplectic,
 then $T_{x}P =T_{x}Q\oplus \pi^{\#}(T_{x}Q^{\perp})$,
$\forall x\in Q$ \cite{W83}.  It is
simple to see that $V_{Q}=\pi^{\#} (TQ^{\perp})$ is a
complementary of $TQ$ in $T_{Q}P$ which possesses all the
required  properties for $Q$ being a \Q submanifold.

Conversely, let $Q$ be a submanifold which satisfies the conditions (i)-(ii) as in
Proposition \ref{pro:sufficient} with the property that
 $\pi ' \in \gm (\wedge^{2} V_{Q})$ is non-degenerate.
 For any $x\in Q$,  it is clear that $\pi^{\#}(T_{x}Q^{\perp})
=\pi'^{\#}(T_{x}Q^{\perp})\subset V_{x}$. Since $\pi '$ is non-degenerate,
$\pi'^{\#}: T_{x}Q^{\perp}\lon V_{x}$ is an isomorphism.
Thus we have $V_{x}=
\pi^{\#}(T_{x}Q^{\perp})$.  Thus it follows that 
$T_{x}P =T_{x}Q\oplus \pi^{\#}(T_{x}Q^{\perp})$. Hence $Q$ 
is cosymplectic. \qed

\begin{cor}
\label{cor:cosymplectic}
Cosymplectic submanifolds are   \Q submanifolds.
\end{cor}
\pf According to Lemma \ref{lem:co}, it  suffices
to verify the last condition (iii) in Proposition \ref{pro:sufficient}.
Since  $V_{Q}=\pi^{\#} (TQ^{\perp})$,
$\gm (V_{Q}) $ is spanned by the vector fields $gX_{f}|_Q$ where $f, g
\in C^{\infty}(P)$ and $f$ is constant along $Q$.
Now clearly $\pr_* [gX_{f} , \pi ]=\pr_* (X_{f}\wedge X_{g}) =0$,
and therefore the condition (iii) in Proposition \ref{pro:sufficient}
is satisfied.  This concludes the proof.\qed       


The following proposition gives a nice characterization for a \Q submanifold.  

\begin{pro}
\label{pro:functions}
Assume that there is a set of functions $f_1, \cdots , f_{k} \in C^{\infty}(P)$
which defines a  coordinate system  on $Q$. Then $Q$ is a \Q submanifold
if
\begin{enumerate}
\item the Hamiltonian vector field $X_{f_{i}}, \ \forall i$, is tangent to $Q$;
\item $d\{f_i, f_j \}\cong 0 \ \ \ (\mbox{mod}\ df_i$) along $Q$.
\end{enumerate}
\end{pro}
\pf Let $V_{Q}=\{v\in T_{Q}P|vf_{i}=0, \  \forall i=1, \cdots , k\}$.
 Clearly, $V_{Q}$
is a vector bundle such that $T_{Q}P=TQ\oplus V_{Q}$.
Moreover $V_{Q}^{\perp}=\mbox{span}\{df_{i}|_{Q}, \ i=1, \cdots , k\}$.
Thus from (i) it follows that $  \pi^{\#}(V_{Q}^{\perp} )\subset TQ$.    
Combining with (ii), we see that $V_{Q}^{\perp}$ is indeed
 a Lie subalgebroid of $T^* P$. Thus $Q$ is a \Q submanifold.\qed

\begin{rmk}
{\em It is natural to ask what is  the rule of the subbundle $V_{Q}$
in the definition of  a \Q submanifold. For a given \Q submanifold,
is $V_{Q}$ unique? If not, what is the relation between
different choices of $V_{Q}$? 

Let $Q$  be a \Q submanifold and 
$\pi |_{Q}=\pi_{Q}+\pi '$, where
$\pi_{Q}\in \gm (\wedge^{2}TQ )$ and $\pi '\in \gm (\wedge^{2} V_{Q})$.
 Assume that
there is another decomposition $T_{Q}P=TQ\oplus V_{Q}'$
satisfying the condition of Definition \ref{def:Q}.
 Then $V_{Q}'$  must correspond to a bundle map $\phi: V_Q \lon TQ$, 
i.e., $V_{Q}' =\{\phi (v) +v |\forall v\in  V_Q\}$.
It is simple to see that the condition (ii') in Remark
\ref{rmk:sufficient} implies that  $\phi \smalcirc (\pi ' )^{\#}
=0$. In paricular, if $Q$ is cosymplectic, $\phi$ must be zero
so $V_{Q}$ is unique. However, in general, it is not clear
how to elaborate  the other condition (iii) in order
to give a clean description of $\phi$.}
\end{rmk}

\subsection{Examples}                        
Now we will discuss some examples of \Q submanifolds.
By Corollary  \ref{cor:cosymplectic}, we already know that
cosymplectic manifolds are \Q submanifolds. The following 
gives a list of other examples.

\begin{ex}
\label{ex:2.13}
{\em Assume that $P$ is a symplectic manifold. If $Q$ is  a 
\Q submanifold, then $Q$ must be a symplectic submanifold
according to Theorem  \ref{thm:Q} (vi). On the other hand,
symplectic submanifolds are automatically \Q  submanifolds since
they are cosymplectic. 	In other words, \Q submanifolds
of  a symplectic manifold are precisely symplectic submanifolds.}
\end{ex}

Another extreme case is the following

\begin{ex}
\label{ex:2.14}
{\em If  $x$ is a point where the Poisson tensor vanishes,
then $\{x\}$ is a \Q submanifold.}
\end{ex}

\begin{ex}
\label{ex:2.15}
{\em Let $P=\reals^n$ be  equipped with a constant Poisson structure.
Then $P$ is a regular Poisson manifold, where  symplectic leaves
are affine subspaces  $x+S$. Here $S$ is the symplectic
leaf through $0$ which is also a linear subspace of  $\reals^n$.
Assume  that an affine subspace $ Q=u+V$ is a \Q
submanifold, where $V$ is a   linear subspace of $\reals^n$. 
By  Theorem \ref{thm:Q} (i),   
we see that $V$ must  admit a complementary subspace $U$ such that
the $P=V\times U$  as  a product of  Poisson  manifolds,
where $V$  and $U$ are  equipped with the constant Poisson structures
$\pi_{Q}(u)$ and $\pi' (u)$ respectively.   This condition is equivalent to
that the intersection of $V$ with  $S$ is a symplectic subspace
of $S$. Conversely, given any such a linear subspace $V$,
then one can decompose  $P=V\times U$ as a product of constant 
Poisson structures. For $Q=V\times \{u\}$, by taking $V_{Q}| \cong Q\times U$ 
to be constant,
one easily sees that the conditions in Proposition \ref{pro:sufficient}
are indeed satisfied. Hence $Q$ is a \Q submanifold.

In conclusion, an affine space $u+V$ is a \Q submanifold iff
$V\cap S$ is a symplectic linear subspace of $S$. }
\end{ex}

The following example, which indicates that being a  \Q
submanifold is indeed a global property,
 was pointed out to the author by Weinstein.

\begin{ex} 
\label{ex:2.16}
{\em Let $P=M\times C$, where each M-slice is a Poisson submanifold.
Namely the Poisson tensor at each point $(x, t)\in M\times C$ is
of the form $\pi (x, t) =\pi_{t} (x)$, where $\pi_{t} (x)$ is a family of
$C$-dependent Poisson structures on $M$. 
Consider  a particular M-slice $Q=M\times \{t_{0}\}$ which is a 
Poisson submanifold. We will investigate when $Q$ becomes a \Q submanifold.

Since we are only  concerned with a small neighborhood of $t_{0} $
in $C$, we may identify $C$ with $\reals^n$ by choosing
a local coordinate system $(t_1 , \cdots , t_n )$.
If $Q$ is a \Q submanifold, then $T_{Q}P=TQ\oplus V_{Q}$ for
some vector bundle $V_{Q}$ along $Q$. Hence $V_{Q}$
must be of the form:
$$V_{Q}=\{  \parrr{}{t_{i}}+X_{i}|i=1, \cdots , n\},$$
where $X_i, \ i=1, \cdots , n, $ are some vector fields on $M$. Clearly 
Condition (ii) in Proposition \ref{pro:sufficient} is satisfied
automatically. Thus according to Remark \ref{rmk:sufficient},
for $Q$ to be a \Q submanifold, it suffices that
$\pr_* [ \parrr{}{t_{i}}+X_{i} , \ \pi_{t} (x)]=0$, for $i=1, \cdots , n$,
which is equivalent to
\begin{equation}
\label{eq:SM}
\parrr{\pi_{t} (x)}{t_{i}}|_{t=t_0} =- [ X_{i}, \pi_{t_0} (x)], \ \ i=1, \cdots , n.
\end{equation} 
This equation precisely means that $\parrr{\pi_{t} (x)}{t_{i}}|_{t=t_0}$ 
is a coboundary with respect to the Poisson cohomology operator
defined by $\pi_{t_0}$.  Thus we conclude that

{\em $Q$ is a \Q submanifold iff the map $f: T_{t_{0}}C\lon H^2_{\pi_{t_0} }(M):
v\lon [v(\pi_{t} )]$ vanishes.}

 Note that $v(\pi_{t} )$ is always a 2-cocycle
with respect to the  Poisson cohomology operator $[\pi_{t_0},  \ \cdot]$
because of the identity $[\pi_{t}, \pi_{t}]=0$. 

As a special case, let us consider the situation where
all $M$-slices are symplectic leaves. Then one obtains 
a map  $\phi : C\lon H^2 (M)$ by taking the   symplectic
class of the fiber. On the other hand, it is known that 
$H^2_{\pi_{t_0} }(M)$ is  canonically isomorphic to $H^2 (M)$.
By identifying these two cohomology groups, we have
\begin{equation}
\label{eq:f}
f=-\phi_{*}.
\end{equation} 
 
To see this relation, let $\omega_t$ denote
the leafwise symplectic forms,  and let
 $\omega_t^{b}: TM \lon T^*M $ and $\pi_{t}^{\#}:  T^*M \lon TM$
be the induced bundle maps by $\omega_t$ and $\pi_{t}$,
respectively.   It follows from the equation
 $\omega_t^{b}\smalcirc \pi_{t}^{\#}=1$ 
that $(v  (\pi_{t}))^{\#}=- \pi_{t}^{\#}\smalcirc (v  (\omega_t  ))^{b}
 \smalcirc \pi_{t}^{\#}$, for any $v\in T_{t_0}C$. Equation  (\ref{eq:f})
thus follows immediately.

Hence  we conclude that a symplectic leaf $M\times\{t_0 \}$ is a 
\Q  submanifold iff $t_0$ is a critical point of the map
$\phi $.
For instance, the symplectic leaves in the Lie-Poisson
$\mathfrak{su}(2)$ can never be \Q submanifolds except for
the zero point.}
\end{ex}

\begin{ex}
\label{ex:2.17}
{\em Let $P=\frakg^*$ be a Lie-Poisson structure corresponding to a Lie
algebra $\frakg$. Consider an  affine space $Q=\mu +  V$.
Assume that $Q$ is  a \Q submanifold where $V_Q$ can be taken
constant. This amounts to saying that we have  a decomposition
$\frakg =\frakl \oplus \frakm$  such that 
$V=\frakm^{\perp}$ and $V_Q\cong Q\times \frakm$ as a vector bundle.
Let  $\{e_{1}, \cdots , e_{k}\}$ be a  basis of $\frakl $
 and  $\{m_{1}, \cdots , m_{t}\}$
a basis of $\frakm$. Then $\{ e_{1}, \cdots , e_{k}, m_{1}, \cdots , m_{t}\}$
consists of  a basis of $\frakg$.  Now let $\{\lambda_{1}, \cdots ,\lambda_{k},
 r_{1},  \cdots , r_t\}$ be its corresponding linear coordinates on $\frakg^*$.
Thus  their Poisson brackets are given by 
$$\{\lambda_i , \lambda_j\}
=a_{ij}^k \lambda_k+b_{ij}^k r_k,\  \{\lambda_i , r_j \}=c_{ij}^k \lambda_k+ d_{ij}^k r_k ,$$
where $a_{ij}^k ,  b_{ij}^k, c_{ij}^k, d_{ij}^k$ are  constants.
It is clear that  $\{\lambda_{1}, \cdots , \lambda_{k}\}$ is a set
of coordinate functions on $Q$ such that $V_{Q}^{\perp}$
is spanned by $d\lambda_{i}, i=1, \cdots , k$.
 Since $d\{\lambda_i , \lambda_j\}
=a_{ij}^k d\lambda_k +b_{ij}^k dr_k$, Condition (ii) 
of Proposition  \ref{pro:functions} implies that $b_{ij}^k =0$.
 On the other hand, we have
$$X_{\lambda_i }|_{Q}=(\{\lambda_i , \lambda_j\}\parrr{}{\lambda_j}
+\{\lambda_i , r_j \}\parrr{}{r_j} )|_{Q}= 
(\{\lambda_i , \lambda_j\} \parrr{}{\lambda_j} 
+(c_{ij}^l \lambda_l +d_{ij}^l \mu_l ) \parrr{}{r_j}  )|_{Q},$$
where $\mu_l =r_{l}(\mu ), \ l=1,  \cdots , t$.
It thus follows that $X_{\lambda_i }$ is tangent to $Q$ iff $c_{ij}^k=0$ 
and $d_{ij}^k \mu_k=0$. The latter  is equivalent to
$\la ad^{*}_{e_j} \mu , m_{j}\ra =0.$ 
Therefore we conclude that $\mu +  \frakm^{\perp}$ is a \Q submanifold
with constant  $V_{Q}$ iff  
$\frakg =\frakl \oplus \frakm$ is a
reductive decomposition (i.e., $\frakl$ is a Lie subalgebra
and $[\frakl ,  \frakm ]\subseteq \frakm $),
 and  $ad^{*}_{\frakl} \mu \in \frakm^{\perp}$. In this case, the
induced Poisson structure can be identified with the Lie-Poisson
structure on $\frakl^*$.}
\end{ex}

\subsection{Local \Q submanifolds}

\begin{defi}
A submanifold $Q$ of a Poisson manifold $P$  is called a local  \Q submanifold
if at each point of $Q$ there is an open neighborhood which is
a \Q submanifold.
\end{defi}

Immediately we have
\begin{pro}
A local  \Q submanifold naturally  carries an induced Poisson structure.
\end{pro}

\begin{ex}
{\em If  $Q$ is  a symplectic leaf of $P$, by Weinstein
splitting theorem \cite{W83},
locally $P\cong Q\times N$ as a product Poisson manifold. It thus follows that
$Q$ is a local \Q submanifold. }
\end{ex}

The following proposition gives a characterization of
local  \Q submanifolds.

\begin{pro}
\label{pro:local}
A submanifold $Q$ of a Poisson manifold $P$ is a \Q submanifold if there
exist local coordinates $(x_1, \cdots , x_{k}, y_{1}, \cdots , y_{t})$
of $P$ at any point $q\in Q$ such that $Q$ is defined by 
$ y_{1}=\cdots =y_{t}=0$ and the Poisson bracket between  coordinate
functions satisfy:
\begin{equation}
\lambda_{ij}(x, 0)=0, \  \forall 1\leq i \leq k, \ 1\leq j \leq t;
 \ \ \parrr{\phi_{ij}}{y_{l}}(x, 0)=0,
 \ \ \forall 1\leq i, \  j  \leq k, \ 1\leq k \leq t,
\end{equation}
where $\phi_{ij}=\{x_i, x_{j}\}$, $\forall 1\leq i, \  j  \leq k$,
 and $\lambda_{ij}=\{x_{i}, y_{j}\}$, $\forall 1\leq i \leq k, \ 1\leq j \leq t$.
\end{pro} 
\pf  Assume that $Q$ is a local  \Q submanifold. Given any point $q\in Q$,
there exists   an open neighborhood $U$ of
$q$ in  $P$   such that $U\cap Q$ is a \Q  submanifold. Let $V_{U\cap Q}$ 
denote the  subbundle as in  the decomposition (\ref{eq:decom}).
By shrinking $U$ to a smaller open neighborhood of $q$ if necessary, one
may always choose  local  coordinates    
$(x_1, \cdots , x_{k}, y_{1}, \cdots , y_{t})$ of $U$ such 
that $U\cap Q$ is defined by
$ y_{1}=\cdots =y_{t}=0$  and $V_{U\cap Q}$ is spanned by
$\{\parrr{}{y_i }|i=1, \cdots , t\}$.
In other words, $\{x_{1}, \cdots , x_{k}\}$ is a set of
coordinates on $Q$ such that $V_{U\cap Q}^{\perp}$ is spanned by
$\{dx_{i}|i=1, \cdots , k\}$.
Then we have
$$d\{x_i, x_{j} \}|_{Q}=\parrr{\phi_{ij}}{y_{l}}(x, 0)dy_{l},  \ (\mbox{mod } dx_{i});\ \ 
X_{x_{i}}|_{Q}=\lambda_{ij}(x, 0)\parrr{}{y_j},  \ (\mbox{mod } \parrr{}{x_{i}}).$$
It thus follows that $\lambda_{ij}(x, 0)=0, \   1\leq i \leq k, \ 1\leq j \leq t$
and $ \parrr{\phi_{ij}}{y_{l}}(x, 0)=0,   1\leq i, j \leq k, \ 1\leq l \leq t$.

Conversely, if   such local coordinates exist in an open neighborhood $U$ of 
$q$ in  $P$, one can show that $U\cap Q$ is a \Q  submanifold.\qed

The following result reveals a connection between local \Q subamnifolds
and  transverse Poisson structures \cite{W83}.

\begin{pro}
If $Q$ is a local \Q submanifold which is a cross section of a
symplectic leaf $S$  at a point $q$
(i.e., $Q$ has complementary dimension to $S$
and intersects with  $S$ at a single point $q$ transversely),
then the induced Poisson structure on $Q$  in a neighborhood of $q$ is
isomorphic to the transverse Poisson structure. 

Conversely, if $Q$ is a cross section of a
symplectic leaf $S$  at a point $q$, then $Q$ 
is a  \Q subamnifold in a  neighborhood of $q$ and 
the induced Poisson structure is isomorphic to the 
transverse Poisson structure.
\end{pro}  
\pf From Weinstein splitting theorem  \cite{W83}, it  follows that a
cross section of a symplectic leaf $S$  must be
a \Q submanifold in a small neighborhood of the intersection point.
It remains to show that the induced Poisson structure
on $Q$ as  a \Q subamnifold  is indeed isomorphic to the transverse
Poisson structure. 

We choose local coordinates as in the proof of  Proposition \ref{pro:local}.
Thus $X_{x_i}$ are all tangent to $Q$ for $i=1, \cdots , k$.
By definition, the transverse Poisson structure
is $\{x_i , x_j\}|_{Q}=\phi_{ij}(x, 0)$, which is precisely
the induced Poisson structure on $Q$ as a \Q submanifold. \qed

An immediate consequence,  by combing with  Example \ref{ex:2.17},
is the following theorem of Molino \cite{M}
and Weinstein \cite{W83}.

\begin{cor}
Let $\mu \in \frakg^*$ and $\frakg_{\mu}$ be  the isotropic Lie algebra 
at $\mu$. If $ \frakg$ admits a reductive decomposition:
 $ \frakg =\frakg_{\mu} \oplus \frakm_{\mu}$, then the transverse Poisson
structure at $\mu$ to the symplectic leaf $G\cdot \mu$ 
(i.e., the coadjoint orbit through $\mu$) is
isomorphic to the Lie-Poisson structure on $\frakg_{\mu}^*$.
\end{cor}

\section{Properties of \Q submanifolds}

This section is devoted to the further study on  properties of \Q submanifolds.

\subsection{Relative modular vector fields}

First we want to see how modular class of a \Q submanifold is related
to that of $P$.  We start with the following: 
 
\begin{lem}
\label{cor:flow}
Let  $Q$ be a \Q submanifold of a Poisson manifold $P$.
Assume that $f\in C^{\infty}(P)$ satisfies the property
 $df|_{Q}\in V_{Q}^{\perp}$.
Denote by $\phi_{t}$ the flow generated by the Hamiltonian vector
field $X_{f} $. Then 
both $TQ$ and $V_{Q}$  (hence $TQ^{\perp}$ and
$V_{Q}^{\perp}$) are stable  under $\phi_{t}$.
\end{lem}
\pf It is clear that $X_{f}$ is tangent to $Q$, and therefore
$[X_{f}, Y]|_{Q}$ is well-defined for any $Y\in \gm (T_{Q}P)$.
If $Y\in \gm (TQ)$, clearly $[X_{f}, Y]|_{Q}\in \gm (TQ)$. Hence 
$TQ$ is stable under $\phi_{t}$.

 Now assume that $Y\in  \gm ( V_{Q} )$.
Let $\tilde{Y}\in \calx (P)$ be any of its extension.
By the graded Jacobi identity, we have
$$[X_{f}, \tilde{Y}]=[[\pi , f], \tilde{Y}]
=[[f, \tilde{Y}], \pi ]-[[\tilde{Y}, \pi ], 
f].$$ 
Now $[[f, \tilde{Y}], \pi ]=-[\tilde{Y}(f), \pi ]=-[\tilde{Y}(f) , 
\pi '' + \tilde{\pi }' ]=-[\tilde{Y}(f) , \pi '']-[\tilde{Y}(f) , 
\tilde{\pi }']$, where
$\pi ''$ and $\tilde{\pi }'$ are  bivector fields on $P$
 as in the proof of Lemma \ref{lem:piQ}, i.e., $\pi '' |_{Q}
=\pi_Q$ and $\tilde{\pi }'\in \gm (\wedge^2 V_{Q} )$.
 Since $\tilde{Y} (f)|_{Q}=
Y(f )|_{Q}=0$, it is obvious that $[\tilde{Y}(f) , \pi '' ]|_{Q}=
[\tilde{Y}(f) , \pi_{Q}]=0$. Thus $[[f, \tilde{Y}], \pi ]|_{Q}
\in \gm (V_Q )$.
On the other hand, according to Equation (\ref{eq:pr}), we have
$\pr_* [\tilde{Y}, \pi ] =[\pr_* \tilde{Y}, \pi_{Q}]=0$. Therefore,
one can write $[\tilde{Y}, \pi ]|_{Q}=\sum Z_{i}\wedge Z_{i}'$ with
$Z_{i}'\in \gm ( V_{Q} )$. Then $[[\tilde{Y}, \pi ], f]|_{Q}=\sum 
(Z_{i}(f)  Z_{i}' -Z_{i}'(f)Z_{i})=  \sum  
Z_{i}(f)  Z_{i}'\in \gm ( V_{Q} )$, since $Z_{i}'(f)=0$ by assumption.
This shows that $[X_{f}, Y]|_{Q } =[X_{f}, \tilde{Y}]|_{Q }
\in \gm ( V_{Q} )$, which
implies that $V_Q$ is stable under the flow $\phi_t$.
This concludes the proof. \qed

We are now ready to introduce  the relative modular class.
Let $\Omega' \in \gm (\wedge^{top}TQ^{\perp}) \cong \gm (\wedge^{top} V_{Q}^* ) $
be a nonzero section, which we always assume exist.
Otherwise, one needs to consider densities as in \cite{W:modular}.
For any $f \in C^{\infty} (Q)$, let $\tilde{f} \in  C^{\infty} (P)$  be
an extension of $f$ satisfying the property $d\tilde{f}|_{Q}  \in V_{Q}^{\perp}$.
Then according to Lemma \ref{cor:flow}   
 the Hamiltonian flow of $X_{\tilde{f}}$ preserves 
both vector bundles $TQ$ and $V_{Q}$, hence it preserves
 $TQ^{\perp}$. It thus follows that $L_{X_f} \Omega'$ 
is a section of $\wedge^{top}TQ^{\perp}$, and therefore
$(L_{X_{\tilde{f}}} \Omega')/\Omega'$ is a  well-defined function on $Q$.
Also it is clear that this function only depends  on the 1-jet
of $\tilde{f}$ along $Q$, and therefore is independent of
the extension. Thus one obtains a  linear map 
$\nu_r :  C^{\infty} (Q) \lon C^{\infty} (Q)$,
$f\lon (L_{X_{\tilde{f}}} \Omega')/\Omega' |_{Q}$. From the
fact that $X_{\tilde{f}}|_Q$ is tangent to $Q$ and
$\Omega' \in \gm (\wedge^{top}TQ^{\perp})$,  one  can easily
show that 
$$\nu_r (fg )=f \nu_r (g) +g \nu_r (f), \ \forall f, g \in C^{\infty} (Q). $$
Hence $\nu_r $ is a vector field on $Q$, which will be called
{\em the relative modular vector field} corresponding to $\Omega'$.

\begin{pro}
\label{pro:rel}
$\nu_r$ is a Poisson vector field with respect to $\pi_Q$.
For different choices of $\Omega' $, the corresponding
relative modular vector fields
$\nu_r$ differ by a Hamiltonian vector field.
\end{pro}  

As a consequence,  $[\nu_r ]$ is a well defined class
 in the Poisson cohomology $H_{\pi_Q}^{1}(Q)$,
which will be called {\em the relative modular
class} of the \Q submanifold $Q$. The proof of   Proposition \ref{pro:rel}
follows from the  lemma below.

Choose a nonzero section  $\Omega_{Q}
\in \gm  (\wedge^{top}V_{Q}^{\perp}) \cong \gm (\wedge^{top}T^* Q)$,
which we again assume exist.
Then $\Omega =\Omega_{Q} \wedge \Omega' \in \gm( \wedge^{top} T^* P|_{Q})$
is a nonzero section.
Extend $\Omega$ to a volume form on $P$ (at least locally along
the submanifold $Q$), which will be denoted by the same symbol
$\Omega$. By $\nu_P$ and $\nu_Q$, we denote the modular
vector fields of the Poisson manifolds $P$ and  $Q$
corresponding to $\Omega$ and $\Omega_{Q}$, respectively.

\begin{lem}
The modular vector fields are related by
\begin{equation}
\label{eq:modular}
\nu_{r}= \pr_* \nu_P -\nu_{Q}.
\end{equation}
\end{lem}  
\pf $ \forall f \in C^{\infty} (Q)$, let $\tilde{f} \in  C^{\infty} (P)$  be
an extension of $f$ satisfying the property $d\tilde{f}|_{Q}  \in V_{Q}^{\perp}$.
Then $L_{X_{\tilde{f}}} \Omega|_{Q}=\nu_P ( \tilde{f} ) \Omega |_{Q}=
(\pr_* \nu_P ) (f)\Omega|_{Q}$, and  $L_{X_{\tilde{f}}}\Omega_{Q} |_{Q}=
\nu_Q (f)\Omega_{Q}$.
>From the derivation law: $L_{X_{\tilde{f}}} \Omega
=(L_{X_{\tilde{f}}}\Omega_{Q})\wedge \Omega'
+\Omega_{Q}\wedge L_{X_{\tilde{f}}} \Omega'$, it follows that
$\pr_* \nu_P (f)=\nu_{Q}(f)+\nu_{r}(f)$.  Equation (\ref{eq:modular})
thus follows.\qed

Another consequence, besides  Proposition \ref{pro:rel}, is
the following:

\begin{pro}
The modular classes of the Poisson structures on $P$ and $Q$ are related by
\begin{equation}
\pr_* [\nu_P ] -[\nu_{Q}] =[\nu_{r}],
\end{equation}
where $\pr_* : H^{1}_{\pi}(P)\lon  H^{1}_{\pi_Q}(Q)$ is the 
morphism as in  Corollary \ref{cor:2.6}.
\end{pro}

\begin{rmk}
{\em It would be interesting to see  how other 
characteristic  classes \cite{C, F} on  $P$ and $Q$
are related, and in particular, how to describe
$\pr_* [C_{k}(P) ] -[C_{k}(Q)]\in H^{*}_{\pi_Q}(Q)$
for other characteristic class $C_{k}$.}
\end{rmk}

\subsection{Poisson actions}
Next we  consider  Poisson actions on \Q submanifolds.
As we shall see below,  \Q submanifolds indeed
 behave nicely under   Poisson group actions, which include
the usual Hamiltonian actions as a special case.

\begin{thm}
\label{thm:Poisson-act}
Assume that $(P, \pi)$ is a Poisson manifold which  admits a  Poisson action of
a  Poisson group $G$.  Assume that $Q$ is a  \Q submanifold stable under the
$G$-action. Then the action of $G$ on $Q$ is also a Poisson action. Moreover,
if $J: P\lon G^*$ is a momentum map, then $J|_{Q}: Q\lon G^*$
is a  momentum map of the $G$-action on $Q$.
\end{thm}
\pf Let $\mu_{P}: T^*P \lon \frakg^*$ and
$\mu_{Q}: T^*Q \lon \frakg^*$  be the linear morphisms dual to the
infinitesimal $\frakg$-actions on $P$ and $Q$, respectively.
Since the infinitesimal $\frakg$-action on $Q$:
 $\frakg \lon \calx (Q)$
is the composition of the infinitesimal $\frakg$-action on $P$:
 $\frakg \lon \calx (P)$  with the projection $\pr_* :\calx (P)
\lon \calx (Q)$, it follows that $\mu_{Q}=\mu_{P} \smalcirc \pr^* $,
where $\pr^* : T^*Q \lon T^* P$ is the dual of the
projection $\pr :T_{Q}P\lon TQ$.
Since $\pr^*$ is  a Lie algebroid morphism according to
Theorem \ref{thm:Q} (iii),  it
follows immediately from Proposition 6.1 in \cite{Xu}
that the $G$-action on $Q$ is also a Poisson action.

Assume that $J: P\lon G^* $ is a momentum map for the Poisson $G$-action
 \cite{Lu}.
I.e., for any $\xi \in \frakg$, $\pi^{\#}(J^* \xi^{l})=\hat{\xi}$,
where $\xi^{l}\in \Omega^{1}(G^* )$ is the left invariant one-form
corresponding to $\xi$, and $\hat{\xi}\in \calx (P)$ is the
vector field on $P$ generated by $\xi$.
Then we have  $\pr_* \pi^{\#}(J^* \xi^{l}) =\hat{\xi}$
since $\hat{\xi}$ is tangent to $Q$.
On the other hand, it is clear that $\pr_* \pi^{\#}(J^* \xi^{l})=
\pi_{Q}^{\#}(J^* \xi^{l})$. This shows that $J|_Q: Q\lon G^*$ is
indeed a momentum map for the  Poisson $G$-action on $Q$. \qed

\subsection{Symplectic subgroupids}

Finally we consider symplectic groupoids of \Q submanifolds.
As we see below, \Q submanifolds are indeed infinitesimal
version of symplectic subgroupids.

\begin{thm}
\label{thm:gpoid}
 If $\poidd{\gm' }{Q}$ is a symplectic subgroupoid of
a symplectic groupoid $\poiddd{\gm}{P}{}$,
  then $Q$ is a \Q submanifold of $P$.
Conversely, if $P$ is an integrable Poisson manifold with symplectic
groupoid $\gm$ and $Q$ is a \Q submanifold  whose
 corresponding cotangent Lie algebroid $T^* Q$  integrates  
to  a Lie subgroupoid $\gm'$ of $\gm$, then $\gm'$ is a symplectic subgroupoid.
\end{thm}
\pf Assume that $\poidd{\gm' }{Q}$ is a symplectic subgroupoid of
a symplectic groupoid $\poiddd{\gm}{P}{}$.
 By $\omega$ and $\omega'$ we denote the symplectic
forms on $\gm$ and $\gm'$ respectively, and by
$A$ and $A'$, we  denote their corresponding  Lie algebroids.
Then $A'$ is a Lie subalgebroid of $A$.
As vector bundles, $A\cong T_{P}^{\alpha}\gm$ and
$A'\cong T_{Q}^{\alpha}\gm' $, and the Lie algebroid
morphism $A'\lon A$ is simply the inclusion: 
$T_{Q}^{\alpha}\gm' \lon T_{P}^{\alpha}\gm$. 
It is well known that $\omega^{b}: T_{P}^{\alpha}\gm \lon T^* P$
  and $(\omega ') ^{b}: T_{Q}^{\alpha}\gm' \lon T^* Q$ are isomorphisms
 of Lie algebroids,  where  $T^* P$ and $T^* Q$ are equipped with
the cotangent Lie algebroids corresponding to the induced
Poisson structures. Thus, one obtains a Lie algebroid
morphism $\phi: T^* Q \lon T^*P$ so that the following diagram
\begin{equation}                         \label{eq:gm}
\matrix{&& {} &&\cr
        &  T_{Q}^{\alpha}\gm' &\vlra& T_{P}^{\alpha}\gm&\cr
        &&&&\cr
   (\omega ') ^{b} &\Bigg\downarrow&&\Bigg\downarrow& \omega^{b}\cr
        &&&&\cr
        &T^* Q&\vlra& T^* P&\cr
        && \phi  &&\cr}
\end{equation}
commutes. In particular, $  \phi (T^* Q )$ is a
 Lie subalgebroid of $T^*P$. In what follows, we 
will  show that $\phi^* \smalcirc i$ is the identity map,
 where $i: T Q\lon TP$ is the inclusion.

Let $\xi \in T_{x}^{*}Q$ be any covector.  Assume that $\xi =(\omega ') ^{b}u$
for some $u\in T_{q}^{\alpha} \gm' $. Then using the commuting diagram (\ref{eq:gm}),
 we have,  for any $v\in T_{x}Q$,
$$\la (i^* \smalcirc  \phi ) \xi , v\ra =\la \phi \xi , v\ra 
=\la  \phi (\omega ') ^{b}u, v \ra
=\la \omega^{b}  u, v\ra  =\omega  (u, v) =\omega' (u, v)
=\la (\omega ') ^{b}u, v\ra =\la \xi, v\ra  .$$
 Therefore $i^* \smalcirc  \phi =id$, or
equivalently $\phi^* \smalcirc i =id$. Let $V_{Q}=ker \phi^*$, 
which is a subbundle of $T_Q P$.
Then $T_Q P=T Q \oplus V_{Q}$.
In fact $V_{Q}^{\perp}=\phi (T^* Q)$, so $V_{Q}^{\perp}$ is
a Lie subalgebroid of $T^*P$. Hence $Q$ is a \Q submanifold.

Conversely, assume that $Q$ is a \Q submanifold
of $P$,  and $\phi =\pr^*: T^*Q\lon T^*P$ is the Lie
algebroid morphism as in Theorem \ref{thm:Q} (iii).
 Let $\gm '\subset \gm $ be a Lie subgroupoid integrating
the Lie subalgebroid $\phi (T^*Q )$. For any $x\in Q$,
we have $T_{x} \gm =T_{x}P\oplus T_{x}^{\alpha }\gm $
and $T_{x} \gm' =T_{x}Q\oplus T_{x}^{\alpha }\gm' $.  By identifying
$T_{x}^{\alpha }\gm $ with $T^*_x P$ via $\omega^{b}$ as
above, one obtains a  decomposition $T_{x} \gm \cong T_{x}P\oplus T^*_x P$, under
which the symplectic form $\omega_{x}\in \wedge^2 T^*_x \gm$ takes
the form:
\begin{equation}
\pmatrix{
0&  I\cr 
-I & \pi (x)\cr}
\end{equation}
Now $T_{x}P=T_x Q \oplus V_{x}$ and $T_{x}^{*}P=
T_x Q^{\perp}\oplus  V_{x}^{\perp} \cong V_{x}^{*}\oplus
T_x^* Q$.  Thus 
$T_{x} \gm \cong  T_x Q \oplus V_{x} \oplus V_{x}^{*}
\oplus T_x^* Q$.
It is clear that  under this decomposition
$T_{x} \gm'  $ corresponds to the subspace
$T_x Q \oplus T_{x}^* Q$. Thus the restriction of $\omega_{x}$
to the subspace $T_{x} \gm' $ has  the form:
\begin{equation}
\pmatrix{
 0 & I \cr -I &\pi_{Q}(x) \cr },
\end{equation}
which is clearly non-degenerate.
It follows immediately that  the pull back of the symplectic
form $\omega $ is non-degenerate  along the identity section $Q$.
To show its non-degeneracy at every point of $\gm'$, it suffices
to show that through each point of $\gm'$, there  exists
a  Lagrangian (local)  bisection $S$ of $\gm$ such that
$S|_Q$ is a bisection of $\gm'$. This is true since any closed
one-form on $Q$ extends to a closed one-form on $P$.\qed

\section{Poisson involutions}

This section is devoted to the study on a special
class of \Q submanifolds arising  as the stable locus of a 
Poisson involution. 
In particular, we discuss Poisson involutions on  Poisson groupoids
as well as  on  Poisson groups. As we will see, such  involutions do
often exist.
Examples include
the standard Poisson group structures on semi-simple Lie groups, Bruhat
Poisson structures on compact semi-simple Lie groups, and
Poisson groupoids connecting with dynamical $r$-matrices of semi-simple
Lie algebras.

\subsection{Stable locus  of  a Poisson involution}

Recall that a Poisson involution on a Poisson manifold $P$
is a Poisson diffeomorphism $\Phi : P\lon P$ such that $\Phi^2 =id$.
Another important class of \Q manifolds arises as follows.

\begin{pro}
Let $\Phi : P\lon P$ be a Poisson involution.
Then its stable  locus $Q$ is a \Q submanifold.
\end{pro}
\pf It is well known that $Q$ is a smooth manifold.
 For any $x\in Q$,  since the linear morphism $\Phi_* :T_x P\lon T_x P$ 
is an involution,   its eigenvalues are either $+1$ or
$-1$. Let  $V_{x}$ denote the  $-1$-eigenspace
of $\Phi_*$,
and $ V_{Q}=\cup_{x\in Q}V_{x}$. Clearly, $T_x Q$ coincides with
 the  $+1$-eigenspace
of $\Phi_*$,   and $T_{x}P=T_x Q \oplus V_{x}$. Since $\Phi_{*} \pi =\pi$, it is
clear that $\pi|_{Q}=\pi_{Q} +\pi '$, where
$\pi_{Q}\in \gm (\wedge^{2}TQ )$ and $\pi '\in \gm (\wedge^{2} V_{Q})$.
It remains to verify  Condition (iii) in Proposition
 \ref{pro:sufficient}.
For this, notice that any vector field $X$  on $P$ can
be decomposed as $X=X^{+}+X^{-}$ where $\Phi_{*}X^{+}=X^{+}$ and
$\Phi_{*}X^{-}=-X^{-}$. Indeed, 
\begin{equation}
\label{eq:X}
X^{+} =\half (X+\Phi_* X)  \ \ \mbox{ and }
X^{-} =\half (X-\Phi_* X). 
\end{equation}
It thus suffices to prove
that  $\pr_* [X^{-} , \pi ] =0$. This is obvious since
$\Phi_{*} [X^{-} , \pi ] =[\Phi_{*} X^{-} , \Phi_{*} \pi ]=-
[X^{-} , \pi ]$. \qed

\begin{rmk} 
{\em The fact that the  stable  locus of a Poisson involution
inherits a Poisson structure was  already hidden in the work of
Bondal \cite{Bondal2} and Boalch \cite{Boalch} in their study of
the    Poisson structures on Stokes matrices. 
On the other hand, an algebraic version of this fact  appeared in
the work of Fernades-Vanhaecke \cite{FV}.}
\end{rmk}

 The Poisson structure on $Q$ indeed can be described more
explicitly in this case.

\begin{pro}
Let $Q$ be the stable locus of a Poisson involution
$\Phi : P\lon P$. Assume that the Poisson tensor $\pi$
on $P$ is $\pi =\sum_{i} X_{i}\wedge Y_{i}$, where $X_{i}$ and $Y_{i}$
are vector fields on $P$. Then the Poisson tensor $\pi_{Q}$ 
on $Q$ is given by $\pi_{Q} =\sum_{i} X_{i}^{+}\wedge Y_{i}^{+} |_{Q}$, where
$X_{i}^{+}$ and $ Y_{i}^{+}$ are defined by Equation (\ref{eq:X}).
\end{pro}

As a consequence of Theorem \ref{thm:gpoid}, we have the following

\begin{cor}
\label{cor:inv-poid}
If $Q$ is the  stable locus of a 
Poisson involution on an  integrable Poisson manifold $P$,
    then $Q$ is always an  integrable Poisson manifold itself.
\end{cor}
\pf Assume that $Q$ is the stable locus of a Poisson 
involution $\Phi: P\lon P$.
 Let $\gm$ be an $\alpha$-connected and simply connected
 symplectic groupoid of  $P$.
To the Poisson involution $\Phi : P\lon P$, there
corresponds to an involutive  symplectic groupoid automorphism
$\tilde{\Phi}: \gm \lon \gm$. Then the stable locus of $\tilde{\Phi}$, which
is a smooth manifold, is a symplectic subgroupoid of $\gm$ 
integrating $Q$. \qed

\subsection{Poisson involutions on Poisson groupoids}

For Poisson groupoids, there is an effective way of producing
Poisson involutions. This is via the so called \sy Poisson groupoids.
Symmetric Poisson groups and their infinitesimal
version: \sy Lie bialgebras, were studied by Fernandes  \cite{Fern, Fern1}
\footnote{Note that our definition here, however,
  is precisely  the opposite to that in \cite{Fern, Fern1}.
We require  that $\Phi$  be group(oid) anti-morphism and Poisson, while
in \cite{Fern, Fern1} $\Phi$ is required to be group morphism and anti-Poisson.}

\begin{defi}
\label{def:sym}
\begin{enumerate}
\item A  \sy Poisson  groupoid consists of a pair $(\gm , \Phi )$,
where $\gm$ is a Poisson groupoid and $\Phi :\gm \lon \gm $ is
a groupiod  anti-morphism which is also a  Poisson involution.
\item A \sy  Lie bialgebroid  consists of a triple  $(A, A^* , \phi)$,
where $(A, A^*)$ is a  Lie bialgebroid and $\phi: A\lon A$
is  an involutive  Lie algebroid anti-morphism
such that $\phi^* : A^* \lon A^*$ is a Lie algebroid morphism.
\end{enumerate}
\end{defi}

\begin{thm}
\label{thm:involution}
Under the assumption that the relevant Lie algebroid is
integrable, there is one-one correspondence between $\alpha$-simply connected
\sy Poisson  groupoids and  \sy  Lie bialgebroids.
\end{thm}
\pf Assume that $(A, A^* , \phi )$ is a \sy  Lie bialgebroid.
Let $\gm $ be an $\alpha$-simply connected Poisson  groupoid corresponding
to the Lie bialgebroid $(A, A^* )$.
 It is known that any Lie algebroid isomorphism
integrates to a Lie groupoid isomorphism for $\alpha$-simply connected
Lie groupoids.
Hence  the Lie algebroid involution
$\phi'=-\phi :  A\lon A$ integrates to a Lie groupoid involution  $\Phi' :
\gm \lon \gm $.
By assumption, $(\phi')^{*}=(-\phi)^* =-\phi^*$ is a Lie algebroid
 anti-morphism. By  the Poisson groupoid duality  \cite{MX1, MX2},
$\Phi'$ is  an  anti-Poisson map.
 Let $\tau : \gm \lon \gm  $ be
 the  map: $\tau (g)=g^{-1}, \ \forall g\in \gm $, which
is clearly a groupoid anti-morphism and an anti-Poisson map.
Set $\Phi =\Phi' \smalcirc \tau $.
Then  $\Phi $ is an integration of $\phi$, which possesses
all the required properties.  

Conversely, if $(\gm , \Phi )$ is  a \sy Poisson
groupoid, then it is clear that  $(A, A^* , \phi)$
is  a \sy Lie bialgebroid, where $\phi: A\lon A$
is the derivative of $\Phi$.  \qed

\begin{rmk}
{\em Note that the roles of $A$ and $A^*$ can be switched for
a \sy Lie bialgebroid. Namely, if $(A, A^* , \phi)$
is  a \sy Lie bialgebroid, then $(A^* , A, -\phi^* )$
is also a \sy Lie bialgebroid. This means that from
a \sy Lie bialgebroid one can in fact  construct a 
pair of Poisson involutions: one on $\gm$ and 
the other on its dual Poisson groupoid  $\gm^*$ (provided
that both $A$ and $A^*$ are integrable).}
\end{rmk}

Theorem  \ref{thm:involution} indicates that a useful source of
producing Poisson involutions on Poisson groupoids
is to construct \sy Lie bialgebroids.
Next we will consider the case of coboundary Lie bialgebroids \cite{LiuX:95},
namely those Lie bialgebroids $(A, A^* )$ where the 
Lie algebroid on the dual $A^*$ is generated by 
an r-matrix  $\Lambda\in \gm (\wedge^2 A)$  with
the property $[X, [\Lambda , \Lambda ]]=0, \ \forall X\in \gm (A)$.

\begin{pro}
 A coboundary Lie bialgebroid  $(A, A^* )$  with an $r$-matrix 
$\Lambda\in \gm (\wedge^2 A)$ is a \sy Lie bialgebroid if there is 
an involutive  Lie algebroid anti-morphism $\phi: A\lon A$ such that
$\phi \Lambda =-\Lambda $.
\end{pro} 
\pf Let $d_{*}: \gm (\wedge^* A)\lon \gm (\wedge^{*+1} A)$ be
the exterior differential induced from the Lie algebroid structure
on $A^*$. Then  for any $X\in \gm (\wedge^* A) $,
$d_{*} X=[\Lambda , X]$.  Hence $(\phi \smalcirc d_{*}) X
=\phi [\Lambda , X]=-[\phi \Lambda , \phi X]=
[\Lambda , \phi X]=  (d_{*} \smalcirc \phi ) X$,
which implies that  $\phi \smalcirc d_{*}= d_{*} \smalcirc \phi$.
Hence $\phi^* : A^*\lon A^*$ is a Lie algebroid morphism.\qed
 
\subsection{Symmetric Courant algebroids}

A nice way of understanding  a Lie bialgebroid $(A, A^* )$
is via its double  $E=A\oplus A^*$, which is  a Courant
algebroid \cite{LWX0}. Roughly, a Courant algebroid
is a vector bundle $E\to M$ equipped with  a non-degenerate
symmetric bilinear form  $(\cdot, \cdot )$ of signature $(n, n)$ on the
fibers,  a bundle map $\rho : E\lon TM$, and a bracket $ [\cdot, \cdot ]$
on $\gm (E)$, which satisfy some complicated compatibility
conditions resembling that of a Lie algebroid up to a homotopy.
Lie bialgebroids  precisely correspond to  splittable
Courant algebroids, namely those which admit two
transversal Dirac structures.  We refer the reader to  \cite{LWX0} for details.

\begin{defi}
\begin{enumerate}
\item A \sy Courant algebroid is a Courant algebroid
 $(E, (\cdot, \cdot ), \rho, [\cdot, \cdot ])$
together with an involutive anti-morphism $\chi: E\lon E$, i.e.,
$$\rho \smalcirc \chi =-f_{*} \smalcirc \chi; \ \ (\chi e_{1}, \chi e_{2}) =-(e_{1}, e_{2});
\ \  \mbox{and } \chi[e_{1}, e_{2}]=-[\chi e_{1}, \chi e_{2}] $$ 
for any $e_{1}, e_{2}\in \gm (E)$, where $f: M\lon M$ is the base map
corresponding to $\chi$;
\item A \sy splittable  Courant algebroid is  a \sy Courant algebroid
$(E, \chi )$, such that
 $E$ admits a   pair of $\chi$-stable transversal Dirac structures.
\end{enumerate}      
\end{defi}

\begin{thm}
\label{thm:sy-courant}
There is a  one-one correspondence between \sy Lie bialgebroids
and \sy splittable  Courant algebroids.
\end{thm} 
\pf  Assume that $(A, A^* , \phi )$   is a \sy Lie bialgebroid.
Let $M$ denote  the base of the Lie bialgebroid $(A, A^* )$,
and $a,  a_*$  their anchors respectively. Denote, by $f: M\lon M$,
the  involution on the base manifold corresponding to
$\phi$. Then $\phi^*$ is a bundle map over the same base
map $f: M\lon M$ since $f$ is an involution.
   Let $E=A\oplus A^*$ be the double of the Lie bialgebroid,
which  is a Courant algebroid \cite{LWX0} over
the base manifold $M$, with anchor $\rho=a +a_{*}$.
Define  $\chi: E\lon E$  by
\begin{equation}
\label{eq:chi}
\chi (X  +\xi )=\phi X  -\phi^* \xi ,  \ \ \ \forall X\in A|_{m} \ \mbox{  and }
\xi \in A^* |_{m}.
\end{equation} 
 Then $\chi $ is clearly an involutive bundle map over the base map $f: M\lon M$.
It is also simple to check that $\chi$ anti-commutes
with the anchor on $E$,  and $(\chi e_{1}, \chi e_{2})=-( e_{1},  e_{2})$
for any $e_{1}, e_{2} \in \gm (E)$. It remains to check that
$[\chi e_{1}, \chi e_{2}]=-[ e_{1},  e_{2}]$ 
for any $e_{1}, e_{2} \in \gm (E)$.  To   this end, it suffices to
show that $[\chi X, \chi  \xi ]=-[X , \xi ]$ for any $X\in
\gm (A)$ and $\xi \in \gm (A^* )$.  First, we will need
the following identities:
\begin{eqnarray}
L_{\phi^{*} \xi} \phi X&=&\phi (L_{\xi }X) ; \label{eq:L}\\
L_{\phi X}\phi^{*}\xi &=&-\phi^{*} ( L_{X} \xi ) \label{eq:L1}.
\end{eqnarray}

Note that for any $\eta \in \gm (A^* )$,
\be
&&\la L_{\phi^{*} \xi} \phi X , \eta \ra\\
&=&(a_{*}\phi^{*} \xi )\la \phi X , \eta \ra -\la  \phi X , [\phi^{*} \xi , \eta ]\ra\\
&=&f_{*} (a_{*} \xi )\la \phi X , \eta >-\la  \phi X , [\phi^{*} \xi , \eta ]\ra\\
&=&(a_{*} \xi )\la X, \phi^* \eta \ra -\la X,  \phi^* [\phi^{*} \xi , \eta ]\ra \\
&=&(a_{*} \xi )\la X, \phi^* \eta \ra -\la X, [\xi , \phi^*  \eta ]\ra \\
&=&\la L_{\xi}X , \phi^* \eta \ra \\
&=&\la \phi (L_{\xi}X ) , \eta \ra .
\ee

Equation (\ref{eq:L}) thus follows. Equation (\ref{eq:L1}) 
can be proved similarly.
Now
\be
&&[\chi  X, \  \chi  \xi]\\
&=&-[\phi X , \ \phi^{*} \xi ]\\
&=&L_{\phi^{*} \xi} \phi X -\half d_{*}\la \phi^{*} \xi , \phi X \ra 
-L_{ \phi X}\phi^{*} \xi+\half d\la \phi^{*} \xi , \phi X \ra \ \ \ 
(\mbox{by Equations (\ref{eq:L})-(\ref{eq:L1})})\\
&=&\phi ( L_{\xi}X) -\half \phi  d_{*}\la \xi , X \ra
+\phi^* ( L_{X }\xi )  -\half \phi^*  d \la \xi , X \ra.
\ee

On the other hand,  
\be
&&\chi [X , \xi ]\\
&=&\chi  [(-L_{\xi}X+\half d_{*}\la \xi ,X\ra )+(L_{X}\xi -\half d\la \xi, X\ra )]\\
&=&-\phi (L_{\xi}X ) +\half \phi d_{*}\la \xi ,X\ra -\phi^* (L_{X}\xi ) +
\half \phi^* d \la \xi , X\ra.
\ee
Thus $ [\chi  X, \chi  \xi]=-\chi [X , \xi ]$.

Conversely, assume that $E$ is a splittable Courant algebroid
such that $E=A\oplus A^*$ for a Lie bialgebroid $(A, A^* )$,
and $\chi :E\lon E$ is an involutive anti-morphism preserving
both  components $A$ and $A^*$. Let $\phi=\chi |_{A}: A \lon A$ and
$\psi=\chi |_{A^*}: A^* \lon A^*$. Then both 
$\phi$ and $\psi$ are  involutive Lie algebroid anti-morphisms.
For any $X\in \gm (A)$ and $\xi\in \gm (A^* ) $,
since $(\chi \xi , \chi X )=-(\xi , X)$,  and 
$\chi X = \phi X$,  $\chi \xi =\psi \xi$,
 it follows   immediately that $\phi^* \psi =-id$,  which
implies that $\psi =-\phi^*$. This concludes the proof. \qed

\subsection{Poisson involutions on dynamical Poisson groupoids}

As a special case, we will consider dynamical Poisson groupoids
introduced by Etingof-Varchenko \cite{EV}.
Recall that  a dynamical $r$-matrix is a function $r: \frakh^{*}
\lon \wedge^{2}\frakg$ satisfying:
\begin{enumerate}
\item $r: \frakh^{*} \lon \wedge^{2} \frakg$  is $H$-equivariant;
\item $ \sum_{i} h_{i} \wedge \frac{\partial r}{ \partial\lambda^{i}}
 -\half [r, r]$ is a constant
$(\wedge^{2} \frakg)^{\frakg}$-valued function  over $\frakh^{*}$,
\end{enumerate}          
where  $\frakh\subset \frakg$ is
a Lie subalgebra,  $H$ is the Lie subgroup of $G$ with Lie algebra $\frakh$,
$\{h_{1}, \cdots , h_{k}\}$ is a basis of $\frakh$,
and $\{\lambda_{1}, \cdots , \lambda_{k}\}$ is  its
induced coordinates on $\frakh^*$.    

It is known  \cite{BK-S, LX1} that a dynamical $r$-matrix naturally defines a
coboundary Lie bialgebroid $(A, A^* , \Lambda )$, where
$A=T\frakh^{*} \times \frakg$, and $\Lambda =
\pi_{\frakh^{*}} +\sum_{i=1}^{k} (\frac{\partial}{\partial \lambda_{i} }
\wedge  h_{i})+r(\lambda )\in \gm (\wedge^{2}A)$. Here
$\pi_{\frakh^*}$ is the Lie-Poisson tensor on $\frakh^*$.

The following theorem can be  verified directly.

\begin{thm}
\label{thm:dy-inv}
Let  $r:  \frakh^{*} \lon \wedge^{2} \frakg$ 
be   a dynamical $r$-matrix.  Assume that
$s: \frakg \lon \frakg$  is  an involutive
Lie algebra anti-morphism,  which preserves   $\frakh$ and satisfies
the property  $s ( r(\lambda  ))=-r (s_{\frakh}^* \lambda )$, 
$\forall \lambda \in \frakh^* $.  Here $s_{\frakh}: 
\frakh \lon \frakh$ is the restriction of $s$ to
$\frakh$. Then
 $(T\frakh^{*} \times \frakg,  T^{*}\frakh^{*} \times \frakg^* , \phi )$,
where $\phi=(-Ts_{\frakh}^{*}, s ):
 T\frakh^{*} \times \frakg\lon T\frakh^{*} \times \frakg$,
 is  a \sy Lie bialgebroid.
\end{thm}

\begin{cor}
Under the same hypothesis as in Theorem \ref{thm:dy-inv}, let
$S: G\lon G$ be the group anti-morphism corresponding
to $s$. Then,
\begin{enumerate}
\item $\Phi : \frakh^{*} \times \frakh^{*} \times G\lon 
\frakh^{*} \times \frakh^{*} \times G$, $\Phi  (u, v, g)
= (s_{\frakh}^{*}v, s_{\frakh}^{*} u, S(  g))$, $\forall u, \ v\in \frakh^*$
and $g\in G$, is  a Poisson involution of the 
dynamical Poisson groupoid
$\frakh^{*} \times \frakh^{*} \times G$.
\item $\gm^0 = \{(u, s_{\frakh}^{*} u, g)|\forall u \in \frakh^*, \ 
g\in G^{0}\}$ is a \Q submanifold, where $G^{0} \subset G$ is
 the stable locus of $S$.   
\end{enumerate}
\end{cor}


\begin{ex}
\label{ex:rmatrix}
{\em  Let $\frakg$ be a semi-simple Lie algebra
over $\complex$ of rank $k$ with a Cartan subalgebra $\frakh$.
Let   $\{e_{\alpha}, f_{\alpha},
h_{i}|\alpha\in \Delta_{+},  1\leq i\leq k\}$ 
be a Chevalley basis.
Then 
\[
r(\lambda) \, =\sum_{\alpha \in \Delta_{+}}
d_{\alpha} \coth ({1 \over 2} < \alpha, \lambda >) e_{\alpha} 
\wedge  f_{\alpha} \]
is a dynamical $r$-matrix over $\frakh^*$, 
where $(e_{\alpha}, f_{\alpha})=d_{\alpha}$,  and
 $\coth (x)  =  {e^x + e^{-x} \over e^x - e^{-x}}$
is the hyperbolic cotangent function \cite{EV}.

Let $s:  \frakg \lon \frakg$ be    a $\complex$-linear morphism,
 which,  on generators, is defined as follows\footnote{Note that $-s$
is  precisely the Cartan involution on the split real form.}:
\begin{equation}
\label{eq:rho}
 s  e_{\alpha} =f_{\alpha}, \ s  f_{\alpha}=e_{\alpha},\ s  h_{i} =h_{i}.
\end{equation}
It is clear that $s$ is an involutive  Lie algebra   anti-morphism
and $s|_{\frakh }=id$. Moreover,  it is  also simple
to see that  $s ( r(\lambda )  )=-r(\lambda )$
for any $\lambda \in \frakh^*$. 
Therefore,  according to Theorem \ref{thm:dy-inv},
   $(T\frakh^{*} \times \frakg,  T^{*}\frakh^{*} \times \frakg^* , \phi )$
 is  a \sy Lie bialgebroid, where $\phi :
  T\frakh^{*} \times \frakg\lon T\frakh^{*} \times \frakg$ is given by
$\phi (v, X)=(-v, s X )$, $\forall (v, X)\in T\frakh^{*} \times \frakg$.
Thus one obtains a pair of Poisson involutions on their
corresponding Poisson groupoids $\Phi: \gm \lon \gm$
and $\Psi: \gm^* \lon \gm^* $.
Now $\gm =\frakh^* \times \frakh^* \times G$, and
$\Phi (u, v, g)= (v, u, S  g)$, $\forall u, \ v\in \frakh^*$
and $g\in G$. Hence,  the stable locus
of $S$ is diffeomorphic to $\frakh^* \times G^{0}$, where 
$G^{0}$ is the stable locus of $S$.  It would be interesting to
compute explicitly the induced Poisson structure on $\frakh^* \times
G^{0}$. On the other hand,  it is
quite mysterious what the 
 stable locus of $\Psi$ should look like, since it is even not clear how to
describe the groupoid $\gm^*$.}
\end{ex}

Let $\frakl$ be a reductive  Lie subalgebra of $\frakg$ containing 
$\frakh$,  i.e.,
\begin{equation}
\label{eq:red}
\frakl =\frakh \oplus \oplus_{\alpha \in \Delta'_{+}}
(\frakg_{\alpha }\oplus \frakg_{-\alpha }),
\end{equation}
where $\Delta'_{+}$ is some  subset of  $\Delta_{+}$.

The claim in Example \ref{ex:rmatrix}  in fact holds in a more general
 situation when $\frakh$ is replaced by $\frakl$.

\begin{pro}
Let $\frakl $ be a   reductive  Lie subalgebra of a semi-simple
Lie algebra $\frakg$ as in Equation (\ref{eq:red}), 
and $r: \frakl^*  \lon \wedge^2 \frakg$
 a dynamical $r$-matrix. Then the map $s : \frakg \lon \frakg$
defined by Equation (\ref{eq:rho}) satisfies  the
conditions as in Theorem \ref{thm:dy-inv}, and therefore
$(T\frakl^{*} \times \frakg,  T^{*}\frakl^{*} \times \frakg^* , \phi )$
 is  a \sy Lie bialgebroid. Here $\phi=(-Ts_{\frakl}^{*}, s ):
 T\frakl^{*} \times \frakg\lon T\frakl^{*} \times \frakg$.
\end{pro}
\pf We prove  this proposition by using the classification
result in  \cite{EV}.
Let $r_{0} : \frakh^* \lon \wedge^2 \frakg$ be  the function:
$$r_{0}  (\lambda )=\sum_{\alpha \in \Delta'_{+}} \frac{1}{(\alpha , \lambda )}
e_{\alpha}\wedge f_{\alpha} .$$

According to \cite{EV}, $\tilde{r}=r|_{\frakh^*}+ r_{0}  : \frakh^* \lon
\wedge^2 \frakg$ is a classical dynamical $r$-matrix on $\frakh^*$.
Hence, from Example \ref{ex:rmatrix}  (the rational case can also
 be  similarly checked),
we know that
 $s (\tilde{r}(\lambda ))=-\tilde{r}(\lambda )$, $\forall \lambda\in
\frakh^*$,
which in turn implies that $s (r(\lambda ))=- r(\lambda ), \ \forall \lambda 
\in \frakh^*$.

Now assume that $\mu=Ad_{x^{-1}}^* \lambda \in \frakl^*$, where
$\lambda \in \frakh^*$ and $x\in L$. Then
\be
s (r(\mu ))&=& s  (r( Ad_{x^{-1}}^* \lambda )) \ \ \mbox{(since $r$ is
$L$-equivariant)}\\
&=& s [Ad_{x} r(\lambda ) ] \\
&=& Ad_{S x^{-1}}s (r( \lambda ))\\
&=& Ad_{S x^{-1}}(-r( \lambda ))\\
&=&-r(Ad_{S x}^{*}\lambda )\\
&=&-r (s^* Ad_{x^{-1}}^* s^* \lambda )\\
&=&-r (s^* Ad_{x^{-1}}^*  \lambda )\\
&=&-r ( s^* \mu ).
\ee
Here we used the identities:
 $s \smalcirc Ad_{x} = Ad_{S x^{-1}} \smalcirc s$
and $Ad_{S x}^{*} =s^* Ad_{x^{-1}}^* s^* $.
Since those  points  $\mu =Ad_{x^{-1}}\lambda$, $\forall 
\lambda \in \frakh^*, \ x\in L$, consist of a dense subset of
$\frakl^*$, the conclusion thus follows immediately. \qed

\section{Poisson involutions on Poisson groups}

In this section,  we turn our attention to Poisson involutions on
 Poisson groups.

\subsection{Symmetric Poisson groups} 
As a special case of Definition \ref{def:sym}, we have 
\begin{defi}
\begin{enumerate}
\item A  \sy Poisson  group consists of a pair $(G , \Phi )$,
where $G$ is a Poisson group and $\Phi :G \lon G$ is
a group  anti-morphism which is also a  Poisson involution.
\item A \sy  Lie bialgebra  consists of a triple  $(\frakg , \frakg^* , \phi)$,
where $(\frakg , \frakg^*)$ is a  Lie bialgebra and $\phi: \frakg \lon \frakg $
is  an involutive  Lie algebra anti-morphism
such that $\phi^* : \frakg ^* \lon \frakg^*$ is a Lie algebra  morphism.
\end{enumerate}
\end{defi}

In this case, a combination of Theorems \ref{thm:involution}
and  \ref{thm:sy-courant} leads to  the following:

\begin{thm}
\label{cor:bialgebra}
\begin{enumerate}
\item   There is a  one-one correspondence between simply connected
\sy Poisson  groups and  \sy  Lie bialgebras.
\item There is one-one correspondence between \sy  Lie bialgebras
 $(\frakg , \frakg^*, \phi )$  and involutive anti-morphisms
$\chi : \sigma \lon \sigma$ (i.e., $(\chi e_{1}, \chi e_{2}) =-(e_{1}, e_{2});
\ \  \mbox{and } \chi[e_{1}, e_{2}]=-[\chi e_{1}, \chi e_{2}] $)
 of the double $\sigma =\frakg\oplus \frakg^*$
 preserving both components $\frakg $ and $\frakg^*$.
\item If $(\frakg, \frakg^*  )$ is  a coboundary Lie bialgebra
with an r-matrix $r\in \wedge^2 \frakg$, then 
$(\frakg, \frakg^* , \phi )$  is a \sy  Lie bialgebra if
$\phi: \frakg \lon  \frakg$ is an involutive  Lie algebra anti-morphism
such that $\phi r =-r $.          
 \end{enumerate} 
\end{thm}

Now  assume that $(\frakg , \frakg^* , \phi  )$ is a \sy Lie bialgebra. 
Then according to the proof of Theorem \ref{thm:sy-courant},
 $\chi: \sigma \lon \sigma, \ \chi (X+\xi )
=\phi X-\phi^* \xi$, $\forall X+\xi \in \frakg\oplus \frakg^*$,
is an involutive  Lie algebra 
anti-morphism, where $\sigma =\frakg\oplus \frakg^* $ is 
the double of the Lie bialgebra.  On the other hand, it is well
known that $(\sigma , \sigma^* )$ itself is a Lie bialgebra with the
$r$-matrix: $r=\sum_{i} X_{i}\wedge \xi^{i} \in \wedge^{2}\sigma$,
 where $\{X_{1}, \cdots , X_{n}\}$ is a basis of $\frakg$ and
$\{\xi^{1},  \cdots , \xi^{n}\}$ is its dual basis of $\frakg^*$.
Then $\chi ( r)=-\sum_{i} \phi X_{i}\wedge \phi^* \xi^{i}  =-r$,
since $\{\phi^* \xi^{1}, \cdots , \phi^* \xi^{n}\}$ is
a dual basis to $\{\phi X_1, \cdots , \phi X_n \}$.
Thus we have proved the following:

\begin{pro}
The double of a  \sy Lie bialgebra is still
a \sy Lie bialgebra.
\end{pro}

\begin{rmk}
{\em  Let $D$  denote the Lie group of $\sigma$.
Then  the same space $D$   possesses three different structures (under
certain  assumptions on completeness): a Poisson group, a symplectic
groupoid  $\gm_G $ over $G$ and a symplectic groupoid $\gm_{G^*}$ over $G^*$.
If $(\frakg, \frakg^* , \phi )$ is a \sy Lie bialgebra,
  then  $\phi$ induces  a Poisson involution on $D$, an involutive automorphism on symplectic
groupoid  $\gm_G $, and an  involutive automorphism on
the symplectic groupoid  $\gm_{G^*}$. These three involutions
are all {\bf different} (see \cite{Bondal2}). Their stable
locuses correspond  to   a \Q submanifold of $D$, a 
symplectic groupoid over the stable locus of $\Phi$,
and a symplectic groupoid over the stable locus of $\Psi$.
Here $\Phi: G\lon G$ and $\Psi: G^* \lon  G^* $ are
the corresponding involutions induced by $\phi$.}
\end{rmk} 

\subsection{Poisson structures on stable locuses}

Below we outline  a scheme to explicitly compute the Poisson 
tensor on the stable locus $Q$ of the
Poisson involution $\Phi$
 for a \sy Poisson group $(G, \Phi )$.
Since $\Phi$ is an involutive
group anti-morphism, we have
\begin{equation}
\label{eq:Ad}
Ad_{\Phi (x)^{-1}}\smalcirc \phi =\phi \smalcirc Ad_{x}: \frakg \lon \frakg, \ \ \ \forall x\in G.
\end{equation}
  
\begin{defi}
\begin{enumerate}
\item A smooth map $\xi : G\lon \wedge^* \frakg$ 
is said to be $\Phi$-equivariant if
\begin{equation}
\label{eq:phi}
\xi (\Phi (x))=Ad_{\Phi (x)} \phi (\xi (x) ),  \ \ \ \forall x\in G;
\end{equation}
\item It is said to be anti-$\Phi$-equivariant if
\begin{equation}
\label{eq:phi-anti}
\xi (\Phi (x))=-Ad_{\Phi (x)} \phi (\xi (x) ),  \ \ \ \forall x\in G. 
\end{equation}
\end{enumerate} 
\end{defi}

Indeed, any  smooth map $\xi : G\lon \wedge^* \frakg$ can
be decomposed as $\xi =\xi^{+}+\xi^{-}$ such that
$\xi^{+}$ is $\Phi$-equivariant and $\xi^{-}$ 
is anti-$\Phi$-equivariant, where
\begin{eqnarray}
\xi^{+} (x)&=&\half [\xi (x)+\phi (Ad_{\Phi (x)^{-1}} \xi (\Phi (x))];
\label{eq:+}\\
\xi^{-} (x)&=&\half [\xi (x)-\phi (Ad_{\Phi (x)^{-1}} \xi (\Phi (x))]. 
\label{eq:-}
\end{eqnarray}
It is simple to see that $\xi : G\lon \wedge^* \frakg$ is $\Phi$-equivariant
(or anti-$\Phi$-equivariant) iff its right translation
 $r_{x*}\xi (x)$ is a $\Phi$-invariant (or anti-$\Phi$-invariant) multi-vector
field on $G$. 


  Let $\delta : \frakg \lon \wedge^2 \frakg$
denote  the cobracket of the Lie bialgebra $(\frakg, \frakg^* )$, which
is also a Lie algebra 1-cocycle, and let $\lambda :G\lon \wedge^2 \frakg$
be  its corresponding Lie group   1-cocycle.  It is well-known that
$\pi (x)=r_{x*} \lambda (x), \ \ \forall x\in G$, is
the Poisson tensor on the Poisson group $G$.
Since $\pi$ is $\Phi$-invariant,
it thus follows that $\lambda : G\lon \wedge^2 \frakg$ is $\Phi$-equivariant.

\begin{pro}
Assume that the  group 1-cocycle
$\lambda : G\lon \wedge^2 \frakg$ is
 $\lambda =\sum_{i}\xi_{i} \wedge \eta_{i}$,
where $\xi_{i} , \ \eta_{i}: G\lon \frakg$. Then 
$\pi_{Q}(x)=\sum_{i} r_{x*} \xi_{i}^{+}(x) \wedge r_{x*} \eta_{i}^{+}(x)|_{Q}$
is the Poisson tensor on $Q$, where $\xi_{i}^{+}$ and $\eta_{i}^{+}$
are defined as in Equations (\ref{eq:+}-\ref{eq:-}). Moreover, the
symplectic leaves  of $Q$ are the intersection of
$Q$ with dressing orbits of $G^*$.
\end{pro}

When $G$ is a coboundary Poisson group, one can write 
$\pi_{Q}$ more explicitly.

\begin{cor}
\label{cor:Q-Poisson}
Under the same hypothesis as in Theorem \ref{cor:bialgebra},  moreover
assume that $G$ is a coboundary Poisson group with
$r$-matrix $r=\sum_i e_{i}\wedge f_{i}\in 
 \wedge^2 \frakg$.
Then  the Poisson tensor on $Q$ is given by
\begin{equation}
\label{eq:Q}
\pi_{Q}=\frac{1}{4}\sum_{i} (\ceV{e_i}+ \Vec{\phi e_i})
\wedge (\ceV{f_i}+ \Vec{\phi f_i})|_{Q} -
\frac{1}{4}\sum_{i} (\Vec{e_i}+ \ceV{\phi e_i})
\wedge (\Vec{f_i}+ \ceV{\phi f_i})|_{Q} ,
\end{equation}
where  $\ceV{e_i}$ and $\Vec{e_i}$ are the left-
and right-invariant vector fields on $G$, respectively, corresponding to
$e_i \in \frakg$; similarly for $\ceV{f_i}$ and $\Vec{f_i}$, etc.

In particular, if $e_{i}, \  f_{i}$ are chosen
such that $\phi e_{i} =e_{i}$ and $\phi f_{i}=-f_{i}$,
then 
\begin{equation} 
\pi_{Q}=\half \sum_i (\ceV{e_i}+\Vec{e_i})
\wedge  (\ceV{f_i}-\Vec{f_i})|_{Q}.
\end{equation}  
\end{cor}
\pf   It is simple to see, by using Equation  (\ref{eq:+}),
that for any $\xi \in \frakg$,
$\xi^{+} (x)= \half (\xi +Ad_{x}(\phi \xi ))$ and
$(Ad_{x}\xi )^{+} (x)= \half (Ad_{x} \xi +\phi \xi )$.
Hence it follows that
$r_{x*}(Ad_{x}\xi )^{+}(x)=\half (\ceV{\xi}+\Vec{\phi \xi})$
and  $r_{x*} \xi^{+}(x)=\half (\Vec{\xi}+\ceV{\phi \xi})$.
It is well known that for a coboundary
Poisson group   $\lambda (x)=\sum_{i}( Ad_{x}e_{i} \wedge
Ad_{x}f_{i} -e_{i} \wedge f_{i})$.
Equation (\ref{eq:Q}) thus follows immediately. \qed

\subsection{Poisson symmetric spaces}

In what follows, we discuss the relation between the
stable locus of the Poisson involution of a \sy Poisson
group and Poisson symmetric spaces.
A Poisson symmetric space is a symmetric space, which is in the mean time
also a Poisson homogeneous space. Poisson symmetric spaces were
studied systematically by  Fernandes in his Ph. D. thesis \cite{Fern, Fern1}, 
to which we refer the reader for details.

Assume that $(G, \Phi )$ is a \sy Poisson group, and 
$Q=\{g|\Phi (g) =g \}$ is the stable locus of $\Phi$.
The following result is standard (c.f. \cite{RS, S}).
For completeness, we outline a proof below.

\begin{pro}
Any connected component of $Q$ is a symmetric space.
\end{pro}
\pf Let $g_0\in Q$ be any fixed point of $\Phi$, and $Q_{g_0}$ the
connected component of $Q$  through  $g_0$.

Consider the  twisted   $G$-action on  (the space) $G$   given by \cite{RS}:
\begin{equation}
\label{eq:Gaction}
g\cdot x =gx\Phi (g), \  \forall g, x \in G.
\end{equation}a
Since $\Phi$ is a  group anti-morphism, this is clearly an action.
Now $\Phi (g\cdot x )=\Phi (gx\Phi (g))=g \Phi (x) \Phi (g)
=g \cdot \Phi (x )$,
so $Q$ is stable under this action. Therefore in particular
$Q_{g_0}$ is stable as well.
Let $Q_{g_0}'$ denote the $G$-orbit through $g_0$. Then
$Q_{g_0}'$ is a homogeneous space $Q_{g_0}' \cong G/H_{g_0}$, where
$H_{g_0}=\{g|g\in G, \ gg_{0}\Phi (g)=g_{0}\}$ is the
isotropic group at $g_0$. Set
\begin{equation}
\Phi_{g_0}: G\lon G, \ \ \Phi_{g_0} (g)=Ad_{g_0}\Phi (g^{-1}), \forall g \in G.
\end{equation}
Then  $\Phi_{g_0}$ is an involutive group homomorphism, since
$$\Phi_{g_0}^2 (g)= \Phi_{g_0} (Ad_{g_0}\Phi (g^{-1}))
=Ad_{\Phi_{g_0}(g_0)} \Phi_{g_0} (\Phi (g^{-1}))  
=Ad_{g_0^{-1}} Ad_{g_0}\Phi (\Phi (g^{-1}))^{-1}=g, \ \ \forall g \in G.$$
It is clear that $H_{g_0}$ is the stable locus of $\Phi_{g_0}$.
Hence $Q_{g_0}'$ is indeed a symmetric space, and its dimension
equals to the dimension of $-1$-eigenspace of $\phi_{g_0}$,
where $\phi_{g_0} =-Ad_{g_0}\smalcirc \phi: \frakg \lon \frakg$
is the Lie algebra involution corresponding to $\Phi_{g_0}$.
On the other hand, the tangent space  $T_{g_0} Q_{g_0}$ 
is spanned by those vectors $v\in T_{g_0} G$ such that
$\Phi_* v=v$. By identifying $T_{g_0} G$ with $\frakg$ by
right translations, $T_{g_0} Q_{g_0}$ can be identified
with the subspace of $\frakg$  consisting of
those elements $X$ satisfying $Ad_{g_0}\smalcirc \phi X=X$,
i.e., the $-1$-eigenspace of $\phi_{g_0}$.
Therefore $Q_{g_0}'$  is a submanifold of $Q_{g_0}$ of the
same dimension, so it must be an open submanifold.
Since it is also closed, they must be identical.
This concludes the proof. \qed

We are now ready to prove the following

\begin{thm}
\label{thm:Poisson-hom}
Let $(G, \Phi )$ be  a \sy Poisson group, and
$Q=\{g|\Phi (g) =g \}$  the stable locus of $\Phi$. 
If  the  Poisson tensor $\pi$ on $G$  vanishes
at a point  $g_0\in Q$, then the connected component
$Q_{g_0}$ is a Poisson symmetric space up to a multiplier of
$2$. In particular, the identity component of $Q$
is a Poisson symmetric space.
\end{thm}
\pf Consider the map 
$$f: G\lon Q_{g_0} , \ \ g\lon  g\cdot g_{0} =g g_{0}\Phi (g), \  \forall g \in G. $$
It suffices to prove that $f$ is a Poisson map, where $Q_{g_0}$ is
equipped with the Poisson tensor $2\pi_{Q}$.

First, it is simple to see that
\begin{equation}
\label{eq:f1}
f_{*} \delta_g =R_{g_{0}\Phi (g)}\delta_g +L_{g g_{0}}\Phi_{*}\delta_g,
\ \ \  \forall \delta_g \in T_g G.
\end{equation} 
On the other hand, we have 
\begin{equation} 
\label{eq:f2}
L_{g g_{0}}\Phi_{*}\delta_g =\Phi_{*}(R_{g_{0}\Phi (g)}\delta_g ).
\end{equation}
 To see this, take a curve $g(t)$ starting at  $g$
with $\parrr{}{t}|_{t=0}g(t)=\delta_g $. Since $\Phi $ is an involutive
anti-morphism, we have $gg_{0}\Phi (g (t))=\Phi (g (t) g_{0}\Phi (g))$.
Equation (\ref{eq:f2}) thus follows by taking the derivative
at $t=0$. Combining Equation (\ref{eq:f1}) with Equation (\ref{eq:f2}),
we are thus lead to
\begin{equation}
\label{eq:f3}
f_{*} \delta_g =2(R_{g_{0}\Phi (g)}\delta_g)^{+}.
\end{equation} 
Now write  $\pi (g) =\sum_{ij} \delta_g^i \wedge \delta_g^j$, where $\delta_g^i$,
$\delta_g^j \in T_g G$. Then we have
$$f_{*} \pi (g) =
4\sum_{ij}(R_{g_{0}\Phi (g)}\delta_g^i )^{+} \wedge (R_{g_{0}\Phi (g)}\delta_g^j)^{+}.
$$

On the other hand, from the multiplicity condition of the Poisson
tensor $\pi (g)$, it follows
that
\be
\pi (g g_{0}\Phi (g) )&=&R_{g_{0}\Phi (g)}\pi (g)+ L_{g}\pi (g_{0}\Phi (g) )\\
&=&R_{g_{0}\Phi (g)}\pi (g)+ L_{gg_0}\pi (\Phi (g))\\
&=&R_{g_{0}\Phi (g)}\pi (g)+ L_{gg_0}\Phi_{*}\pi (g)\\
&=&R_{g_{0}\Phi (g)}\pi (g)+\Phi_{*}(R_{g_{0}\Phi (g)}\pi (g))\\
&=&\sum R_{g_{0}\Phi (g)}\delta_g^i \wedge R_{g_{0}\Phi (g)}\delta_g^j
+\sum \Phi_{*}R_{g_{0}\Phi (g)}\delta_g^i \wedge \Phi_{*}R_{g_{0}\Phi (g)a
}\delta_g^j  .
\ee
Here we used the assumption $\pi (g_{0})=0$ in the 
second equality.  Therefore  we have
$$\pi_{Q}( g g_{0}\Phi (g))=2\sum (R_{g_{0}\Phi (g)}\delta_g^i )^{+} \wedge (R_{g_{0}\Phi (g)}\delta_g^j)^{+}.$$
This concludes the proof. \qed

\begin{rmk}
{\em
\begin{enumerate}
\item Theorem \ref{thm:Poisson-hom} would follow from Theorem 
\ref{thm:Poisson-act}, if the action  defined by
Equation (\ref{eq:Gaction}) were a Poisson action where
the Poisson group  is equipped with the Poisson tensor
$\pi (g)$ while the  space it acts, which is $G$ again,
 is equipped with $2\pi (g)$. However, this is false in general.
So we can see that a Poisson group action
on a Poisson manifold  $P$ may not be Poisson
action , but it can still be Poisson when restricted
to the stable locus $Q$.
\item One drawback of Theorem \ref{thm:Poisson-hom} is that
 the stable locuses do  not seem to produce any new examples of
Poisson manifolds  for \sy Poisson groups
in contraction to what one may have initially expected.
A good point, on the other hand, is  that one might be
 able to quantize these Poisson structures
on stable locuses including the one on
Stokes matrices $U_{+}$ (see Example \ref{ex:5.7}) using quantum
homogeneous spaces.
\item One can construct a symplectic groupoid of a Poisson
symmetric space by means of reduction \cite{Xu:92}. On the other
hand, according to  Corollary \ref{cor:inv-poid}, for 
a stable locus Poisson structure, one can construct a
symplectic groupoid directly via the lifted involution
on the corresponding symplectic groupoid.  It is
interesting to compare these two approaches in our case here.
\end{enumerate}}
\end{rmk}

\subsection{Examples}

We end the paper with a list of examples.    We refer the reader to
\cite{Fern} for a  complete list of orthogonal \sy Lie bialgeras, which
also contains examples below.

\begin{ex}
\label{ex:5.7}
{\em Let $\frakg$ be a semi-simple Lie algebra
of   rank $k$ over $\complex$  with a Cartan subalgebra $\frakh$.
Let   $\{e_{\alpha}, f_{\alpha},
h_{i}|\alpha\in \Delta_{+},  1\leq i\leq k\}$
be a Chevalley basis.  It is well known that $(\frakg , \frakg^* )$
is a coboundary  Lie bialgebra with   $r$-matrix:
\[
r \, =\sum_{\alpha \in \Delta_{+}}
d_{\alpha}  (e_{\alpha} \wedge  f_{\alpha}), \]
where $d_{\alpha} =(e_{\alpha}, f_{\alpha})$.

As in Example \ref{ex:rmatrix}, let $\phi  : \frakg \lon \frakg$ be    the 
$\complex$-linear morphism, which, on generators, is defined as follows: 
$$\phi e_{\alpha} =f_{\alpha}, \ \phi  f_{\alpha}=e_{\alpha},\ \phi  h_{i} 
=h_{i}. $$
It is clear that $\phi $ is an involutive  Lie algebra   anti-morphism
and $\phi r =-r$. Therefore,  $(\frakg , \frakg , \phi )$ is
a \sy Lie bialgebra, which in turn induces a pair of 
\sy Poisson groups $(G, \Phi )$ and $(G^* , \Psi )$.
Thus one obtains a pair of Poisson
involutions: $\Phi : G\lon G$  and $\Psi  :G^* \lon G^*$,
which are the group anti-morphisms corresponding to the 
Lie algebra anti-morphisms: $\phi  : \frakg \lon \frakg$
  and $-\phi^* : \frakg^* \lon \frakg^*$, respectively. 

For $\frakg =\mathfrak{sl}(n, \complex )$,  it is well-known that
$G=SL(n, \complex )$ and $G^* =B_{+}*B_{-}$. It is
simple to see that $\Phi$ and $\Psi$ are given by the following:
$$\Phi: SL(n, \complex )\lon SL(n, \complex ), \ \ \ 
\Phi ( A) =A^T, \ \forall A\in SL(n, \complex );$$
and 
 $$\Psi: B_{+}*B_{-}\lon B_{+}*B_{-},  \ \ \  
\Psi (B, C)=(C^{T}, B^{T}), \  \forall (B, C)\in B_{+}*B_{-}.$$

The stable locus  of $\Phi$ thus  consists of all symmetric matrices in
$SL(n, \complex )$. On the other hand,  the set  $U_{+}$  of Stokes matrices
 (i.e.  upper triangular matrices with all main diagonal entries being $1$)
 can be identified with the identity   component of the  stable locus of $\Psi $.
As a consequence, both the space  $S$ of
symmetric matrices in
$SL(n, \complex )$ and the space  $U_{+}$ of Stokes matrices  admit
natural  Poisson structures. These Poisson manifolds, together
with their symplectic groupoids, were studied
in details by Bondal \cite{Bondal2} in connection with 
his study of triangulated categories. Independently,
the Poisson structure on  $U_{+}$
 was also obtained independently by Dubrovin \cite{Dubrovin} in the
$3\times 3$-case and  Ugaglia
\cite{Ugaglia} in  general  in connection with the study of
Frobenius manifolds. From
a very different aspect,  
the relation between  the Poisson structure on  the space of
Stokes matrices $U_{+}$ and the  Poisson group  $B_{+}*B_{-}$   was independently found 
 by Boalch in his study of the so called ``monodromy map" \cite{Boalch}.
 We refer the reader to \cite{Bondal2, Boalch}
for details.  As a consequence of Theorem \ref{thm:Poisson-hom}, we conclude
that both $S$ and $U_{+}$ are indeed Poisson symmetric spaces.

\begin{thm}
 Up to a multiplier $2$,
\begin{enumerate}
\item  the map $SL(n, \complex )\lon S, \ A\to AA^{T}, 
\forall  A\in SL(n, \complex ) $, is a Poisson map. Indeed $S$ is a Poisson 
symmetric space with  the  Poisson $SL(n, \complex )$-action:
$$SL(n, \complex )\times S\lon S, \ \ A\cdot X = AXA^{T}, \ \forall A
\in SL(n, \complex ), X\in S;   $$
\item  the map $B_{+}*B_{-} \lon U_{+}, \ (B, C)\to BC^{T},
\forall (B, C)\in B_{+}*B_{-}$ is a Poisson map. Indeed $U_{+}$
is a Poisson symmetric space with  the  Poisson $B_{+}*B_{-}$-action:
$$B_{+}*B_{-} \times U_{+}\lon  U_{+},\ \ (B, C)\cdot X =BXC^{T},
(B, C)\in B_{+}*B_{-}, X\in U_{+}.$$
\end{enumerate} 
\end{thm}

}
\end{ex}

\begin{ex}
{\em   Let $K$ be a compact semi-simple Lie group
with Lie algebra $\frakk$, and 
$\frakt$  its  Cartan subalgebra. It is well known that
$K$ admits a standard Poisson group structure called Bruhat
Poisson structure \cite{LW}.
 Let $\frakg =\frakk^{\complex}$ be its complexification, which
is a  complex semi-simple Lie algebra. Choose a  Chevalley basis
$\{e_{\alpha}, f_{\alpha},
h_{i}|\alpha\in \Delta_{+},  1\leq i\leq k\}$ of $\frakg$
 as in Example \ref{ex:5.7},
 then  $\{X_{\alpha}, Y_{\alpha}, t_{i}|
\alpha\in \Delta_{+},  1\leq i\leq k\}$, where 
\begin{equation}
X_{\alpha}=e_{\alpha} -f_{\alpha}, \ Y_{\alpha}= \ii (e_{\alpha}+ f_{\alpha}),
\ \mbox{ and }
t_{i}=\ii h_{i},
\end{equation}
 is a basis (over $\reals$) of $\frakk$, and
\begin{equation} 
\hat{r}=\ii r= \ii\sum_{\alpha \in \Delta_{+}}
d_{\alpha} (  e_{\alpha} \wedge  f_{\alpha})=
\sum_{\alpha \in \Delta_{+}} \half d_{\alpha} 
X_{\alpha}\wedge Y_{\alpha}
 \in \wedge^2 \frakk 
\end{equation}
is indeed
the $r$-matrix generating the corresponding Lie bialgebra
$(\frakk, \frakk^* )$. Let $ \phi: \frakg \lon \frakg$ be the
anti-morphism as in Example \ref{ex:rmatrix}. It is then clear that
$\phi (X_{\alpha})=-X_{\alpha}, \ \phi (Y_{\alpha} )=Y_{\alpha}$,
and $\phi (t_{i})=t_{i}$, so $\frakk$ is stable
under $\phi$. 
It is also clear that $\phi \hat{r} =-\hat{r}$.
Hence $(\frakk, \frakk^* ,  \hat{\phi})$, where 
 $\hat{\phi}=\phi|_{\frakk}: \frakk \lon \frakk$,
is a \sy Lie bialgebra.
 Thus it induces a pair
of Poisson involutions $\hat{\Phi}: K\lon K$ and
$\hat{\Psi}: K^*\lon K^*$. 

To describe the stable locuses of these involutions,
 we need to  consider the double of the Lie bialgebra
$(\frakk, \frakk^* )$, which is  isomorphic to $\frakg$ as a real Lie algebra.
According to Theorem \ref{cor:bialgebra},
  $\hat{\phi}$ induces an involutive Lie algebra
antimorphism (over $\reals$)   $\chi: \frakg\lon \frakg$,
under which both $\frakk$ and $\frakk^*$ are stable and
whose restrictions to these Lie subalgebras are $\hat{\phi}$ and
$-\hat{\phi}^*$, respectively. In our case, a straightforward
computation yields that  on generators $\chi$ is given by:
\be
&&\chi (\ii e_{\alpha})=\ii e_{\alpha}, \  \ \chi (\ii f_{\alpha})=\ii 
f_{\alpha},
\  \ \chi (\ii h_{i})=\ii h_{i}\\
&&\chi (e_{\alpha})=-e_{\alpha}, \  \ \chi (f_{\alpha})=-f_{\alpha},
\  \ \chi (h_{i})=-h_{i}.
\ee
In other words, $\chi (X)=-\bar{X}$, $\forall X\in \frakg$.
On the group level, $\chi$ induces  an involutive Lie group
antimorphism $\Upsilon: G\lon G$ such that $\Upsilon (g)=\bar{g}^{-1}$,
$\forall g\in G$,
where $G$ is a simply connected Lie group (considered as a real
Lie group) integrating the Lie algebra $\frakg$.
By $Q$, we denote the stable locus of $\Upsilon$, i.e.,
$Q=\{g\in G|\bar{g}=g^{-1} \}$. Then the stable
locus of $\hat{\Phi}$ and $\hat{\Psi}$ are $K\cap Q$ and
$K^* \cap Q$, respectively. In particular, according to
Corollary \ref{cor:Q-Poisson}, 
\begin{equation}
\label{eq:piK}
\pi_{Q}=
\sum_{\alpha \in \Delta_{+}} \frac{1}{4} d_{\alpha}
(\Vec{X_{\alpha}}-\ceV{X_{\alpha}})\wedge (\ceV{Y_{\alpha}}+\Vec{Y_{\alpha}})
\end{equation}
is the Poisson tensor on $K\cap Q$. Theorem  \ref{thm:Poisson-hom}
implies that the map $g\to g\bar{g}^{-1}$ is indeed  Poisson maps
(up to a factor of $2$)   when being restricted to $K$ and $K^*$.

For $K=SU(n)$, its dual group $K^*$ is isomorphic to
$SB(n, \complex )$, and the double $G\cong SL(n, \complex )$,
considered as a real Lie group.
Thus $Q= \{A\in SL(n, \complex )|\bar{A}A=I\}$. Hence
we have $K\cap Q\cong \{A |A^{*}A=\bar{A}A=I, \det A=1\}$,
which is the submanifold of $SU(n)$ consisting  of all  symmetric matrices.
On the other hand, $K^* \cap Q\cong \{A\in SB(n, \complex )|\bar{A}A=I\}$.

We note that $SB(n, \complex )$ is Poisson  diffeomorphic
to the linear Poisson structure on $\mathfrak{sb} (n,  \complex )$
according to Ginzburg-Weinstein theorem \cite{GW}. The
recent result of Boalch  \cite{Boalch} suggests that there may exist
a Poisson diffeomorphism $SB(n, \complex )\lon \mathfrak{sb} (n,  \complex )$
commuting with the Poisson involutions, where
the Poisson involution on $SB(n, \complex )$ is
given by $A\lon \bar{A}^{-1}$ while the Poisson involution on 
$ \mathfrak{sb} (n,  \complex )$
 is: $A\lon -\bar{A}$.  If so,
the induced Poisson structures on their stable
locus should be isomorphic.  The latter is a lot easier
to compute and in fact is  a linear Poisson structure.}
\end{ex}

\end{document}